\documentclass[12pt]{article}

\usepackage{amssymb}
\usepackage{amsmath}
\usepackage{amsbsy}
\usepackage{amscd}
\usepackage{amsfonts}
\usepackage{amsthm}
\usepackage{mathrsfs}
\usepackage{hyperref}
\hypersetup{hidelinks}
\usepackage{fullpage}
\usepackage{mathdots}
\usepackage{graphicx,subfigure}
\usepackage[english]{babel}
\usepackage[utf8]{inputenc}
\usepackage{tikz}
\usepackage{adjustbox}
\usepackage{authblk}
\usepackage{caption}
\usepackage{multicol}

\usepackage{mathtext}
\usepackage{stackengine,scalerel}

\usepackage{tikz-cd}
\usetikzlibrary{backgrounds,fit, matrix}
\usetikzlibrary{positioning}
\usetikzlibrary{calc,through,chains}
\usetikzlibrary{arrows,shapes,decorations,automata, petri}
\usepackage[boxsize=1.25em, centerboxes]{ytableau}

\usetikzlibrary{trees}

\usetikzlibrary{arrows}

\newtheorem{Theorem}[equation]{Theorem}
\newtheorem{Corollary}[equation]{Corollary}
\newtheorem{Lemma}[equation]{Lemma}
\newtheorem{Proposition}[equation]{Proposition}

\theoremstyle{definition}

\newtheorem{Example}[equation]{Example}

\newtheorem{Remark}[equation]{Remark}

\numberwithin{equation}{section}
\numberwithin{figure}{section}

\usepackage{mathtools}
\newcommand{\E}{\mathbb{E}}

\newcommand{\1}{\mathbf{1}}
\newcommand{\dG}{d_{X,G,P}}

\newcommand{\Bin}[1]{B_{#1}(G)} 
\newcommand{\indeg}{\deg_{\mathrm{in}}}
\newcommand{\wtG}{\omega_G}
\newcommand{\ballG}{B_G}

\newcommand{\F}{{\mathbb F}}

\newcommand{\Z}{{\mathbb Z}}

\newcommand{\R}{{\mathbb R}}

\usepackage[all,2cell]{xy}
	\usepackage{amssymb, hyperref}
	\usepackage[all,2cell]{xy}
	\usepackage{hyperref}
	\usepackage{array}
	\usepackage[capitalise]{cleveref}

\usepackage{changes}

\begin{document}
	
\title{Mixing and Merging Metric Spaces Using Directed Graphs}

\author[1]{Mahir Bilen Can}
\author[2]{Shantanu Chakrabartty}

\affil[1]{\small{Department of Mathematics, Tulane University, USA\\{\it mcan@tulane.edu}}}
\affil[2]{\small{Department of Electrical and Systems Engineering, Washington University in St. Louis, USA\\{\it shantanu@wustl.edu}}}    

\date{}

\maketitle

\begin{abstract}
Let $(X_1,d_1),\dots,(X_N,d_N)$ be metric spaces with $d_i:X_i\times X_i\to[0,1]$ and let $X=\prod_{i=1}^N X_i$.
{
Let $G$ be a directed graph with vertex set $[N]=\{1,\dots,N\}$ and edge set $E(G)$.
Given edge weights $P=\{p_{ij}\in(0,1]:(i,j)\in E(G)\}$, define the closed in-neighborhood
$B_j(G):=\{j\}\cup\{i:(i,j)\in E(G)\}$ and put $p_{jj}=1$.
} 
Define, for $g,h\in X$,
\[
d_{X,G,P}(g,h)\;:=\;1-\frac{1}{N}\sum_{j=1}^N\ \prod_{i\in B_j(G)}\bigl(1-d_i(g_i,h_i)\bigr)^{1/p_{ij}}.
\]
We prove that $d_{X,G,P}$ is a metric and analyze its behavior under graph operations (disjoint unions, Cartesian products).
Two quantitative principles govern this family: (i) \emph{edge densification} (adding edges) increases the distance pointwise and in expectation, with an exact single-edge increment formula and a uniform $1/N$ sensitivity bound; and (ii) \emph{weight monotonicity}: increasing an edge weight $p_{ij}$ decreases the metric.
{
When the coordinate pairs are independent across coordinates, we obtain a closed-form factorization of the mean distance.}
In the discrete specialization $X=\mathbb{F}_q^N$ with $d_i\in\{0,1\}$ and $G$ acyclic, {
after replacing $G$ by the oriented reflexive transitive closure $G^{\preceq}$ of Corollary~\ref{intro:C1}, the unnormalized distance $D_{X,G^{\preceq},P}:=N d_{X,G^{\preceq},P}$ coincides
} with the size of the hereditary closure of the support, yielding a poset/digraph \emph{sphere-packing bound} via the counting of lower ideals.
{
Finally, we give a continuous analog of the construction, yielding a metric on $C([0,1],\mathbb C)$ through an explicit integral formula.}
These results unify and generalize several graph-based metrics used in coding theory while providing new design levers for mean-distance control.

\medskip


\medskip


\end{abstract}

\section{Introduction}

The main goal of this article is to study a graph-based construction that mixes and merges multiple
metric spaces through a directed graph with edge weights. Throughout, we work with metric spaces denoted by 
$(X_1,d_1),\dots,(X_N,d_N)$, where $d_i:X_i\times X_i\to[0,1]$ is the distance function.
For notational economy, we use $X$ to denote the product set $X_1\times\cdots\times X_N$.
Let $[N]$ denote $\{1,\dots, N\}$. 
{
Given a directed graph $G$ on $[N]$ {
with edge set $E(G)$ and
} positive edge weights $P=\{p_{ij}\in(0,1]:(i,j)\in E(G)\}$, set
$B_j(G):=\{j\}\cup\{i:(i,j)\in E(G)\}$ and put $p_{jj}=1$.} 
We define, for $g,h\in X$,
\begin{equation}
d_{X,G,P}(g,h)\;:=\;1-\frac{1}{N}\sum_{j=1}^N\ \prod_{i\in B_j(G)}\bigl(1-d_i(g_i,h_i)\bigr)^{1/p_{ij}}.
\label{eq:main_metric}
\end{equation}


We first prove that $d_{X,G,P}$ is a metric on $X$. We then analyze how this
distance behaves under natural graph operations (disjoint unions and Cartesian products), and in two
areas of interest: (i) discrete coordinates over finite fields, which {
recover
graph-/poset-based metrics used in error-correcting codes
}, and (ii) dense-graph limits via graphons.

\medskip

We develop (in Section~\ref{S:EEM}) several principles that yield important information about what to look for in designing our metrics for various applications. 
In a nutshell, these principles are 
(i) \emph{edge densification} (adding edges) increases $d_{X,G,P}$ pointwise and in mean, with an
exact single-edge increment identity and a uniform $1/N$ sensitivity bound
(Theorem~\ref{thm:mean-monotone} and Corollary~\ref{cor:edge-sensitivity});
(ii) under coordinate-pair independence, the \emph{mean distance factorizes} into degree-weighted
nodewise summaries $m_i(t)=\E[(1-d_i)^t]$, yielding a closed form and Jensen-type bounds
(Theorem~\ref{thm:mean-factorization} and Corollary~\ref{cor:mean-bounds});
(iii) \emph{weight monotonicity}: increasing an edge weight $p_{ij}$ (in our convention $p_{ij}\in(0,1]$)
\emph{decreases} the metric (Proposition~\ref{P:weight-monotone}). These properties explain and
quantify the right-shift we observe in empirical distance histograms as edges are added.
{
To give an example, let us consider a two-vertex ($N=2$) graph $G$ with the single edge $1\to2$ of weight $w$. Therefore, $P = \bigl(\begin{smallmatrix}1& w\\0&1\end{smallmatrix}\bigr)$, which implies $B_1(G)=\{1\}$ and $B_2(G)=\{1,2\}$.  
Let the metric space be defined as $X= [0,1] \times [0,1]$ with $d_i(u,v)=|u-v|$, for $u,v \in [0,1]$.
For $g,h \in X$ and writing $a=|g_1-h_1|$ and $b=|g_2-h_2|$, we obtain
\[
d_{X,G,P}(g,h)=1-\frac{(1-a)+(1-a)^{1/w}(1-b)}{2}
\]
following the metric definition in equation~\ref{eq:main_metric}.
For $g=(0,0)$ and $h=(1/2,1/2)$, the edgeless graph gives distance $1/2$.
Adding $1\to2$ with $w=1$ gives $5/8$, and decreasing its weight to $w=1/2$ gives $11/16$, both highlighting the \emph{edge densification} and \emph{weight monotonicity} property of the metric space for this simple graph.
}

\smallskip
Now we mention specific corollaries of our results that record sensitivity to the graph structure and to edge weights.

\begin{Corollary}\label{intro:C01}
Let $G-E$ be obtained by removing an edge $E$ from $G$, and $G+E$ by adding a missing edge $E$ to $G$ with an arbitrary weight $p_E \in (0,1]$ added to $P$.
Then for all $g,h\in X$,
\[
d_{X,G-E,P}(g,h)\ \le\ d_{X,G,P}(g,h)\ \le\ d_{X,G+E,P}(g,h).
\]
\end{Corollary}

\medskip

\begin{Corollary}\label{intro:C02}
Fix an existing edge $(i,j)$ and perturb its weight to $p_{ij}\pm\Delta$ with {
$0<\Delta<p_{ij}$ and $p_{ij}+\Delta\leq1$.
}
Then for all $g,h\in X$,
\[
d_{X,G,P+\Delta}(g,h)\ \le\ d_{X,G,P}(g,h)\ \le\ d_{X,G,P-\Delta}(g,h).
\]
\end{Corollary}

These corollaries follow by comparing the terms in the defining average and using $0\le 1-d_i\le1$ with the map
$p\mapsto1-(1-d_i)^{1/p}$ decreasing on $(0,1]$.

\medskip

Discrete specializations of our metrics are useful in coding-theoretic settings. 
{When $X=\F_q^N$ and $d_i\in\{0,1\}$, the unnormalized distance $D_{X,G,P}:=N d_{X,G,P}$ equals the size of the \emph{hereditary closure}, that is the lower order ideal in the poset
induced by $G$, of the support of $g-h$ after replacing $G$ by its oriented reflexive transitive closure. 
\begin{Corollary} \label{intro:C1}
We assume that the underlying sets in the main result, $X_1,\dots, X_N$, are finite and that the corresponding metrics $d_1,\dots, d_N$ take values in $\{0,1\}$. Let $G^{\preceq}$ denote the oriented reflexive transitive closure of $G$, oriented so that the closed in-neighborhood of $j$ is its forward reachability set:
\[
B_j(G^{\preceq})=\{i\in[N]: \text{there is a directed path, possibly trivial, from } j \text{ to } i \text{ in }G\}.
\]
Then, for every $g=(g_1,\dots,g_N) \in X$ and $h =(h_1,\dots, h_N)\in X$, we have 
\begin{align*}
D_{X,G^{\preceq},P}(g,h) = |\langle {\rm supp}(g,h) \rangle|,
\qquad
 d_{X,G^{\preceq},P}(g,h)=\frac{|\langle {\rm supp}(g,h) \rangle|}{N},
\end{align*}
where $D_{X,G^{\preceq},P}:=N d_{X,G^{\preceq},P}$, ${\rm supp} (g,h) = \{ j\in [N] \mid g_j \neq h_j \}$, 
and $\langle {\rm supp}(g,h) \rangle$ is the hereditary closure of the support. 
\end{Corollary}}

This corollary not only yields a clean identification with poset/digraph weights but also {
opens doors
} to many interesting questions in error-correcting codes along the lines of~\cite{BGL, EFM2018}.
We would like to emphasize however that previous approaches do not directly address the combination of various metrics using graphs. 
Our proposed metric $d_{X,G,P}$ and its unnormalized coding-theoretic version $D_{X,G,P}$ (as in Corollary~\ref{intro:C1}, for example) {
allow different normalized coordinate metrics to be combined in one construction. This permits the mathematical treatment of mixed channels and modulation schemes.
} These include phase-shift keying modulation~\cite{Lee1958, Proakis2008}, or channels that are susceptible to synchronization errors~\cite{MBT}.
{
In addition, the classical sphere-packing argument gives a bound expressed via counts of lower ideals, as follows.
}

\begin{Theorem}[Sphere packing in the generalized poset/digraph metric]\label{thm:poset-sphere-packing}
Let $X=\F_q^N$, $d_i\in\{0,1\}$, and suppose $G$ is a DAG. For $S\subseteq[N]$, write $\langle S\rangle$ for the smallest lower (hereditary) set containing $S$ in the poset induced by $G$, and define the $G$-weight $\wtG(x):=|\langle \mathrm{supp}(x)\rangle|$. For $r\in\{0,1,\dots,N\}$ set the (closed) ball
\begin{equation*}
\ballG(r)\;:=\;\{x\in X:\ \wtG(x)\le r\}.
\end{equation*}
If a code $C\subset X$ has minimum unnormalized $G$-distance $D_{\min}$ with respect to $D_{X,G^{\preceq},P}$, then
\begin{equation}
\label{eq:packing-bound}
|C|
\;\le\;\frac{|X|}{\displaystyle\max_{z\in X}\bigl|\,z+\ballG(\lfloor (D_{\min}-1)/2\rfloor)\,\bigr|}
\;\le\;\frac{q^N}{\displaystyle\sum_{s=0}^{\lfloor (D_{\min}-1)/2\rfloor} J_G(s)\,(q-1)^s}\,,
\end{equation}
where $J_G(s)$ denotes the number of lower ideals of cardinality $s$ in the poset of $G$.
\end{Theorem}

{
The first inequality in \eqref{eq:packing-bound} is the usual sphere-packing bound for the poset metric introduced by Brualdi, Graves, and Lawrence~\cite{BGL}.
Translation invariance makes all balls of a fixed radius equally large.
The second inequality uses only vectors whose support is itself a lower ideal, so it can be weaker than the bound using the full ball volume.
For an antichain, $J_G(s)=\binom Ns$, and it reduces to the classical $q$-ary Hamming bound.
}

\medskip

Moving forward, we discuss dense graphon limits. A graphon, also known as a graph limit, is a symmetric measurable function that serves as an analytic representation for sequences of large, dense graphs. Introduced in the foundational work by Lov\'asz and Szegedy~\cite{LovaszBalazs2006}, graphons provide a powerful framework for studying the properties of large graphs and their limiting behavior, analogous to how probability distributions describe the limits of sequences of random variables~\cite{BorgsC17}.

\begin{Corollary} \label{intro:C2}
{Assume that $d:\mathbb C\times\mathbb C\to[0,1)$ is continuous and that $\rho(u,v):=-\log(1-d(u,v))$ is a metric.
If, in the limit as $N \rightarrow \infty$, the graph $G$ from the main result approaches a graphon $W:[0,1]^2 \rightarrow (0,1]$ defined over the random variables $x,y \in [0,1]$, then
\begin{align*}
d_{G,W} (g,h) := 1 - \mathbb{E}_{x}\left[ \exp \left(\mathbb{E}_{y} \left[\frac{\log \left(1- d(g(y),h(y))\right)}{W(x,y)} \right]\right) \right]
\end{align*}
is a metric defined on the set of all continuous functions $g,h:[0,1] \rightarrow \mathbb{C}$. Here, $\mathbb{E}_{x}(\cdot)$ and $\mathbb{E}_{y}(\cdot)$ denote expectations with respect to $x$ and $y$, respectively, whenever the expectations are well defined. }
\end{Corollary}

\medskip

In the rest of the introduction, we proceed with a minor notational change in order not to mix vectors and their entries. 
For aesthetic reasons, we will use $\mathcal{X},\mathcal{G}$, and $\mathcal{P}$ for the product spaces, graphs, and for the edge weights.
The next corollary records the behavior of $d_{\mathcal{X},\mathcal{G},\mathcal{P}}$ under disjoint unions, which are the additive operation in the usual semiring of directed graphs. A {\em semiring} is an algebraic structure $(S,+,\cdot)$, where addition, denoted by $+$, and multiplication, denoted by $\cdot$, are associative binary operations on the set $S$. Furthermore, $(S,+)$ is required to be a commutative monoid with 0, and $(S,\cdot)$ is required to be a monoid with identity 1. Finally, multiplication distributes over addition both on the left and the right. Semirings are one of the most commonly studied structures in algebra. It is not difficult to see that the set of all directed graphs is a semiring, where the addition is defined as the disjoint union of directed graphs, and the multiplication is defined as the Cartesian product of directed graphs.  In this context, our metric has the following additivity property:

\begin{Corollary} \label{intro:C3}
Let $(X_1,d_1),\dots, (X_N,d_N)$ be metric spaces as in the main result. 
Let $\mathcal{G}_1$ (resp. $\mathcal{G}_2$) be a directed graph on the set $\{1,\dots, N_1\}$ (resp. on the set $\{N_1+1,\dots, N\}$). 
We denote by $\mathcal{X}^{(1)}$ (resp. by $\mathcal{X}^{(2)}$) the product metric space $X_1\times \cdots \times X_{N_1}$ (resp. $X_{N_1+1}\times \cdots \times X_N$).
Let $\mathcal{P}_1$ (resp. $\mathcal{P}_2$) be the matrix of probability scores on the edges of $\mathcal{G}_1$ (resp. on the edges of $\mathcal{G}_2$). 
If $\mathcal{G}$ is the disjoint union of subgraphs, $\mathcal{G}= \mathcal{G}_1 \sqcup \mathcal{G}_2$,
then we have 
$$
d_{\mathcal{X},\mathcal{G},\mathcal{P}} \ = \ \frac{N_1}{N}\ d_{\mathcal{X}^{(1)}, \mathcal{G}_1 ,\mathcal{P}_1} \ + \ \frac{N_2}{N}\ d_{\mathcal{X}^{(2)}, \mathcal{G}_2 ,\mathcal{P}_2},
$$
where $\mathcal{X} = \mathcal{X}^{(1)}\times \mathcal{X}^{(2)}$ and $N_2 = N-N_1$. 
\end{Corollary}

A noteworthy special case of our previous corollary is that, if the number of factors in the product-sets $\mathcal{X}^{(1)}$ and $\mathcal{X}^{(2)}$ is the same,  
then the distance  
$$
d_{\mathcal{X},\mathcal{G},\mathcal{P}} ( (\mathbf{g}_1,\mathbf{g}_2), (\mathbf{h}_1,\mathbf{h}_2)),\quad\text{where}\quad
(\mathbf{g}_1,\mathbf{g}_2), (\mathbf{h}_1,\mathbf{h}_2) \in \mathcal{X}^{(1)}\times \mathcal{X}^{(2)},
$$ 
is the average of the distances 
$$
d_{\mathcal{X}^{(1)}, \mathcal{G}_1 ,\mathcal{P}_1}(\mathbf{g}_1,\mathbf{h}_1) \ \text{ and } \  d_{\mathcal{X}^{(2)}, \mathcal{G}_2 ,\mathcal{P}_2} (\mathbf{g}_2,\mathbf{h}_2).
$$ 
In particular, if we have $\mathcal{X}^{(1)} = \mathcal{X}^{(2)}$, $\mathcal{P}_1 = \mathcal{P}_2$,
and $\mathcal{G}_1 = \mathcal{G}_2$, then 
we have 
$$
d_{\mathcal{X},\mathcal{G},\mathcal{P}} ( (\mathbf{a},\mathbf{a}), (\mathbf{b},\mathbf{b})) =
d_{\mathcal{X}^{(1)},\mathcal{G}_1,\mathcal{P}_1} (\mathbf{a},\mathbf{b}) \quad\text{for all}\quad \mathbf{a},\mathbf{b} \in \mathcal{X}^{(1)}.
$$

\medskip
We now proceed to examine the effects of the Cartesian product of graphs on our metrics. 

Let $ \mathcal{G}_1 = (V_1, E_1) $ and $ \mathcal{G}_2 = (V_2, E_2) $ be two finite directed graphs. The {\em Cartesian product} of $ \mathcal{G}_1 $ and $ \mathcal{G}_2 $, denoted by $ \mathcal{G}_1 \square \mathcal{G}_2 $, is a directed graph whose vertex set is given by 
    \[
    V(\mathcal{G}_1\ \square \ \mathcal{G}_2) := V_1 \times V_2.
    \]
The edge set of $\mathcal{G}_1 \square \mathcal{G}_2$, denoted $E(\mathcal{G}_1\square \mathcal{G}_2)$ is defined so that there is a directed edge from $ (u_1, u_2) $ to $ (v_1, v_2) $ in $ \mathcal{G}_1 \square \mathcal{G}_2 $ if and only if 
    \begin{itemize}
        \item $ u_1 = v_1 $ and $ (u_2, v_2) \in E_2 $, or
        \item $ u_2 = v_2 $ and $ (u_1, v_1) \in E_1 $.
    \end{itemize}
In other words, edges in $ \mathcal{G}_1 \square \mathcal{G}_2 $ connect vertices where one coordinate is adjacent in one graph, while the other coordinate remains unchanged.
This dichotomy allows us to associate a well-defined probability score to each edge of $\mathcal{G}_1\square \mathcal{G}_2$ as long as we have probability scores on the edges of the graphs $\mathcal{G}_1$ and $\mathcal{G}_2$ separately. 
Let $p_{u_1,v_1}$ (resp. $p_{u_2,v_2}$) denote the probability score associated with the edge $(u_1,v_1)$ in $\mathcal{G}_1$ (resp. $(u_2,v_2)$ in $\mathcal{G}_2$).
If $((u_1,u_2),(v_1,v_2))$ is an edge of $\mathcal{G}_1\square \mathcal{G}_2$, then we define
\begin{align*}
p_{(u_1,u_2),(v_1,v_2)} := 
\begin{cases}
p_{u_1,v_1} & \text{ if $u_2 = v_2$,}\\
p_{u_2,v_2} & \text{ if $u_1=v_1$.}
\end{cases}
\end{align*}

\bigskip

The last corollary uses $d_{\mathcal{X},\mathcal{G},\mathcal{P}}$ to understand the nature of a product metric space. Let $X_1,\dots, X_N$ be metric spaces. Let $N_1\in \{1,\dots, N-1\}$ and let $N_2:=N-N_1$. For $i\in \{1,\dots, N_1\}$ and $j\in \{N_1+1,\dots, N\}$, we denote by $F(X_i,X_j)$ the set of all functions from $X_i$ to $X_j$. 
The corresponding {\em uniform metric} on $F(X_i,X_j)$ is defined by
\[
d_{i,j}(f, g) := \sup_{x \in X_i} d_{X_j}(f(x), g(x)),
\]
where $ f, g \in F(X_i, X_j) $, and $ \sup $ denote the supremum. We let $F$ denote the product metric space:
$$
F:= \prod_{i=1}^{N_1} \prod_{j=N_1+1}^N F(X_i,X_j). 
$$
It is notationally convenient to view $F$ as a collection of $N_1\times N_2$ matrices 
$$
\mathbf{f} := (f_{ij})_{i=1,\dots, N_1, j=1,\dots, N_2},
$$
where the $(i,j)$-th entry of the matrix $\mathbf{f}$ is a function $f_{ij}$ in $F(X_i,X_{j+N_1})$. 
Our goal is to understand the nature of the joint metric on $F$. 

{
\begin{Corollary} \label{intro:C4}
Maintain the notation from above, and put $M:=N_1N_2=|V(\mathcal{G}_1\square\mathcal{G}_2)|$.
For $\mathbf{g},\mathbf{h}\in F$, set
\[
\delta_{u_1,u_2}(\mathbf{g},\mathbf{h})
:=d_{u_1,N_1+ u_2}(g_{u_1,u_2},h_{u_1,u_2}).
\]
Then the metric associated with the Cartesian product graph is given by the exact formula
\begin{align*}
d_{F, \mathcal{G}_1 \square \mathcal{G}_2, \mathcal{P}}(\mathbf{g}, \mathbf{h})
=1-\frac{1}{M}\sum_{v_1\in V_1}\sum_{v_2\in V_2}
&\left(1-\delta_{v_1,v_2}(\mathbf{g},\mathbf{h})\right)\\
&\times\prod_{(u_1,v_1)\in E_1}\left(1-\delta_{u_1,v_2}(\mathbf{g},\mathbf{h})\right)^{1/p^{(1)}_{u_1v_1}}\\
&\times\prod_{(u_2,v_2)\in E_2}\left(1-\delta_{v_1,u_2}(\mathbf{g},\mathbf{h})\right)^{1/p^{(2)}_{u_2v_2}}.
\end{align*}
\end{Corollary}

Although the formula is genuinely nonfactorizing in general, weighted in-regular Cartesian factors with separable survival data do factor at the level of $1-d$. The following proposition records this positive case.

\begin{Proposition}
\label{P:cartesian-factorization}
We maintain the notation of Corollary~\ref{intro:C4}, and set
\[
s_{u_1,u_2}(\mathbf{g},\mathbf{h}):=1-\delta_{u_1,u_2}(\mathbf{g},\mathbf{h}).
\]
Assume that the survival matrix is multiplicatively separable:
\[
s_{u_1,u_2}(\mathbf{g},\mathbf{h})=r_{u_1}c_{u_2}
\]
for some $r_{u_1},c_{u_2}\in[0,1]$. Define the weighted in-degrees
\[
A(v_1):=\sum_{(u_1,v_1)\in E_1}\frac{1}{p^{(1)}_{u_1v_1}},
\qquad
B(v_2):=\sum_{(u_2,v_2)\in E_2}\frac{1}{p^{(2)}_{u_2v_2}}.
\]
Assume that $\mathcal{G}_1$ and $\mathcal{G}_2$ are weighted in-regular, meaning $A(v_1)\equiv A$ and $B(v_2)\equiv B$ are constant. Define
\[
S_1(r)
:=
\frac{1}{N_1}
\sum_{v_1\in V_1}
r_{v_1}^{1+B}
\prod_{(u_1,v_1)\in E_1}
r_{u_1}^{1/p^{(1)}_{u_1v_1}},
\]
and
\[
S_2(c)
:=
\frac{1}{N_2}
\sum_{v_2\in V_2}
c_{v_2}^{1+A}
\prod_{(u_2,v_2)\in E_2}
c_{u_2}^{1/p^{(2)}_{u_2v_2}}.
\]
Then the Cartesian-product metric satisfies
\[
1-d_{F,\mathcal{G}_1\square\mathcal{G}_2,\mathcal{P}}(\mathbf{g},\mathbf{h})=S_1(r)S_2(c).
\]
Equivalently, if $d_1^{(B)}:=1-S_1(r)$ and $d_2^{(A)}:=1-S_2(c)$, then
\[
d_{F,\mathcal{G}_1\square\mathcal{G}_2,\mathcal{P}}(\mathbf{g},\mathbf{h})
=
d_1^{(B)}+d_2^{(A)}-d_1^{(B)}d_2^{(A)}.
\]
\end{Proposition}

\begin{Remark}
\label{R:in-regular-examples}
When all edge weights are constant, the weighted in-regularity hypothesis of
Proposition~\ref{P:cartesian-factorization} reduces to ordinary in-regularity, that is, the quantity
\[
\deg^{-}(v):=\#\{u:(u,v)\in E\}
\]
is independent of $v$. This condition is satisfied by many standard families of
directed graphs, including directed cycles, regular tournaments, circulant/Cayley
digraphs, complete bidirected graphs, and directed tori with their standard
orientations.
\end{Remark}
}

{
In Proposition~\ref{P:cartesian-factorization}, for a fixed pair $(\mathbf g,\mathbf h)$, separability means that the nonnegative matrix $(s_{ij})$ has rank at most one.
It can be checked by testing $s_{ij}s_{k\ell}=s_{i\ell}s_{kj}$ for every choice of indices.
For example, on a product of two sensor arrays, a model in which
$\delta_{ij}=1-\exp(-a_i-b_j)$ with $a_i,b_j\geq0$ gives $s_{ij}=e^{-a_i}e^{-b_j}$.
This describes separate row and column effects. Interaction terms depending jointly on $i$ and $j$ generally destroy this factorization.
For instance, the survival matrix $\bigl(\begin{smallmatrix}1&1/2\\1/2&1\end{smallmatrix}\bigr)$ is not separable.
Weighted in-regularity is a separate requirement on the two graphs.
The proposition applies when both conditions hold, as in a homogeneous array with separable data.
}

We close our introduction by mentioning areas of interest for applications. 
Our framework is broadly applicable in areas with heterogeneous coordinates and structured
interactions such as machine learning and related data-fusion settings
(see \cite{gulrajani2017, Mohan2023, SuarezDiaz2018}) and ecological systems with networked dependencies
(see \cite{Dove2023, Lazar2020, OzgodeYigin2023}). The quantitative principles that we outlined give actionable levers:
adding edges raises expected distances in a controlled way, while tuning weights lowers or raises distances
predictably, enabling principled architecture choices.

The organization of our work is as follows. 
Section~\ref{S:Preliminaries} reviews directed-graph metrics and graphons. 
Section~\ref{S:Joint} contains the metric proof and preparatory lemmas.
Sections~\ref{S:Specializations}, \ref{S:Graphons}, and \ref{S:Semiring} establish Corollaries~\ref{intro:C1}, \ref{intro:C2}, \ref{intro:C3}, and \ref{intro:C4}, Proposition~\ref{P:cartesian-factorization},
as well as the sphere-packing bound Theorem~\ref{thm:poset-sphere-packing}. 
Section~\ref{S:Examples} presents examples and empirical distributions.
Section~\ref{S:EEM} develops the mean-shift, factorization, sensitivity, and weight-monotonicity results referenced above.
We conclude in Section~\ref{S:Conclusions}.

\section{Preliminaries}\label{S:Preliminaries}

In this section, we introduce our notation, discuss preliminary results, and make comments on the usefulness of our proposed metrics. Throughout this article, $\Z_+$ stands for the semigroup of positive integers.

\subsection{Graphs and Posets.}
{
A \emph{digraph} is a pair $G=(V,E)$, where $V$ is a finite set of vertices and
$E\subseteq V\times V$ is an irreflexive relation,  that is, $(v,v)\notin E$ for
all $v\in V$.  The elements of $E$ are called directed edges.  A \emph{graph}
is a digraph for which $E$ is symmetric, while an \emph{oriented graph} is a
digraph for which $E$ is anti-symmetric.  Throughout the paper, unless otherwise
specified, our directed graphs are digraphs in this sense.  The self-coordinate
which appears later in the closed in-neighborhood $B_j(G)$ is added by convention
and is not regarded as a loop in $G$.

A \emph{directed walk} from $u$ to $v$ is a sequence of vertices
$v_0,v_1,\ldots,v_k$ such that $v_0=u$, $v_k=v$, and
$(v_{i-1},v_i)\in E$ for every $1\leq i\leq k$.  A directed walk is a
\emph{directed path} if the vertices $v_0,\ldots,v_k$ are distinct; we allow the
trivial directed path of length zero from a vertex to itself.  A \emph{directed
cycle} is a sequence $v_0,v_1,\ldots,v_k$ with $k\geq 2$, such that
$v_0=v_k$, the vertices $v_0,\ldots,v_{k-1}$ are distinct, and
$(v_{i-1},v_i)\in E$ for every $1\leq i\leq k$.  A \emph{directed acyclic
graph}, or DAG, is a digraph with no directed cycles.
}

Let $\mathcal{G}_1 = (V_1, E_1)$ and $\mathcal{G}_2 = (V_2, E_2)$ be directed graphs, where $E_1 \subseteq V_1 \times V_1$ and $E_2 \subseteq V_2 \times V_2$. The {\em Cartesian product of directed graphs} $\mathcal{G}_1 \square \mathcal{G}_2$ is a directed graph with the following properties:
\begin{itemize}
    \item The vertex set of $\mathcal{G}_1 \square \mathcal{G}_2$ is the Cartesian product $V_1 \times V_2$.
    \item There is a directed edge from $(u_1, u_2)$ to $(v_1, v_2)$ in $\mathcal{G}_1 \square \mathcal{G}_2$ if and only if one of the following conditions holds:
    \begin{enumerate}
        \item $u_1 = v_1$ and $(u_2, v_2) \in E_2$ (a directed edge in $\mathcal{G}_2$), or
        \item $u_2 = v_2$ and $(u_1, v_1) \in E_1$ (a directed edge in $\mathcal{G}_1$).
    \end{enumerate}
\end{itemize}

A \textit{partially ordered set} (poset) is a pair $P = (S, \preceq)$, where $S$ is a set and $\preceq$ is a binary relation on $S$ that is reflexive ($x \preceq x$ for all $x \in S$), antisymmetric (if $x \preceq y$ and $y \preceq x$, then $x = y$), and transitive (if $x \preceq y$ and $y \preceq z$, then $x \preceq z$).
Similarly to our assumption on the finiteness of graphs, unless otherwise stated, our posets are assumed to be finite.

A \textit{subposet} of $P$ is a subset $T \subseteq S$ together with the restriction of $\preceq$ to $T$.
For elements $x, y \in S$, the \textit{interval} between $x$ and $y$ is defined as $[x, y] = \{z \in S \mid x \preceq z \preceq y\}$. 
If $x \npreceq y$, then $[x, y] = \emptyset$. If $[x, y] = \{x, y\}$, we say that \textit{y covers x}.

A poset $P = (S, \preceq)$ can be conveniently visualized by its \textit{Hasse diagram}, which is a directed graph $(V, E) = (S, C)$, where $C$ is the set of all cover relations.
A DAG $G$ naturally defines a poset. 
If there is a directed edge from vertex $x$ to vertex $y$ in $G$, we write $x \lessdot y$. 
The {
reflexive transitive closure
} of the relation $\lessdot$, denoted by $\leq$, defines a partial order on the vertices of $G$, thus forming a poset.

For a poset $(R, \preceq)$, a \emph{lower order ideal} is a subset $S \subseteq R$ such that for every $r \in R$ and $s \in S$, if $r \preceq s$, then $r \in S$. A closely related concept appears in graph theory.

Let $(V, E)$ be a directed graph. A subset $W \subseteq V$ is called a \emph{hereditary subset} if every vertex $v \in V$ with a path to some vertex in $W$ also belongs to $W$. 

{
When the directed graph is acyclic, hereditary subsets correspond precisely to lower order ideals.
} To see this, consider the poset $(V, \preceq)$ where $v \preceq v'$ if and only if there exists a directed path from $v$ to $v'$ in the graph $(V, E)$. Under this construction, the hereditary subsets of $V$ are exactly the lower order ideals of the associated poset, and vice versa.

Let $A$ be a subset of the vertex set of a directed graph $(V,E)$.
We denote by $\langle A \rangle$ the smallest hereditary subset $W$ of $(V,E)$ such that $A \subseteq W$.

\subsubsection{Metric Spaces.}

A \textit{metric space} is an ordered pair $(X, d)$, where $X$ is a set and $d: X \times X \to [0, \infty)$ is a function called a \textit{metric} (or \textit{distance function}) that satisfies the following axioms for all $x, y, z \in X$:
\begin{enumerate}
    \item \textit{Non-negativity:} $d(x, y) \geq 0$, and $d(x, y) = 0$ if and only if $x = y$. 
    \item \textit{Symmetry:} $d(x, y) = d(y, x)$.
    \item \textit{Triangle inequality:} $d(x, z) \leq d(x, y) + d(y, z)$. 
\end{enumerate}

Let $(X, d_X)$ and $(Y, d_Y)$ be metric spaces, where $X$ and $Y$ are sets and $d_X$ and $d_Y$ are distance functions on $X$ and $Y$, respectively. 
A function $f: X \to Y$ is called an {\em isometry} if for all points $a, b \in X$, the distance between $a$ and $b$ in $X$ is equal to the distance between their images $f(a)$ and $f(b)$ in $Y$. In other words:
\begin{equation}
d_Y(f(a), f(b)) = d_X(a, b).
\end{equation}
Notice that an isometry is necessarily injective. 
Unless otherwise noted, we assume that isometries are surjective as well. 
When we write `$f$ is an isometry of $(X,d)$' we mean that $f: X\to X$ is an isometry. 

A metric $d$ is {\em normalized} if its range is contained within the unit interval: 
$$
0 \leq d(x, y) \leq 1\quad \text{ for all $x, y \in X$}.
$$ 
Note that this is equivalent to $\sup_{x,y \in X} d(x,y) \leq 1$. If $\sup_{x,y \in X} d(x,y) = 1$, we say the metric is \textit{strictly} normalized.

As we mentioned in the introduction, every metric can be normalized to take values in the interval $[0,1]$. 
\begin{Lemma}\label{L:normalized}
Let $d$ be a metric on a set $X$, not necessarily normalized.
If $\sup_{s,t\in X}d(s,t)=0$, set $d'=d$. Otherwise, let $d' : X\times X \to \R$ be the function defined by 
$$
d'(x,y) = 
\begin{cases}
\frac{ d(x,y) }{\sup_{s,t\in X}d(s,t)} & \text{ if $X$ is finite or $d$ is bounded},\\ \\
\frac{ d(x,y) }{1+d(x,y)}& \text{ otherwise}.
\end{cases}
$$
Then $d'$ is a normalized metric. 
\end{Lemma}
\begin{proof}
If $\sup_{s,t\in X}d(s,t)=0$, then $X$ has only one metric class and the zero metric is already normalized. Thus assume $\sup_{s,t\in X}d(s,t)>0$.
If $d'$ is defined by $d'(x,y):=\frac{d(x,y)}{\sup_{s,t\in X}d(s,t)}$ for $x,y\in X$, 
then it is straightforward to check that $d'$ is a normalized metric since it is obtained from $d$ by scaling it by its supremum. 
We proceed with the second case. 
The verification of the following properties are easy to check: 
\begin{enumerate}
\item $d'(x,y)\geq 0$ for every $x,y\in X$, with equality if and only if $x=y$. 
\item $d'(x,y) =d'(y,x)$ for every $x,y\in X$. 
\end{enumerate}
To check the triangle inequality, we apply the following observation:
For three nonnegative real numbers $a,b$, and $c$ such that $a \leq b+c$, we have
\begin{align}\label{A:subadditive1}
\frac{a}{1+a}  \leq \frac{b}{1+b} + \frac{c}{1+c}.
\end{align}
Indeed, it is easy to check that 
$$
a + a(b+c) + abc \leq (b+c+2bc) + a(b+c+2bc)
$$
which is equivalent to the inequality 
$$
a(1+b)(1+c) \leq (1+a)(b(1+c)+c(1+b)).
$$
Dividing the both sides of this inequality by $(1+a)(1+b)(1+c)$ we obtain the inequality in (\ref{A:subadditive1}).
This shows that the triangle inequality holds.
Hence, the proof follows.
\end{proof}

\subsection{Metrics from directed graphs.}

If the underlying set of a metric space $(X, d)$ is a group such that the group operations are compatible with the metric (specifically, translations are isometries), then $(X, d)$ is called a \textit{metric group}.

Let $K$ be a group. The $N$-fold direct product $K^N = K \times \cdots \times K$ admits a natural translation-invariant metric, the \textit{Hamming metric}, defined as follows. Let $e$ denote the identity element of $K$. The \textit{support} of an $N$-tuple $v = (v_1, \dots, v_N) \in K^N$ is defined by
$$
\operatorname{supp}(v) = \{i \in [N] \mid v_i \neq e\}.
$$
The Hamming metric on $K^N$ is then given by
$$
d_H(v, w) = |\operatorname{supp}(vw^{-1})|,
$$
where $vw^{-1} = (v_1w_1^{-1}, \ldots, v_Nw_N^{-1})$. 
It is easy to check that this metric is translation-invariant.
The \textit{Hamming weight} of $v$, denoted by $\omega_H(v)$, is the Hamming distance from $v$ to the identity element $(e, \dots, e) \in K^N$:
$$
\omega_H(v) = d_H(v, (e, \dots, e)) = |\operatorname{supp}(v)|.
$$ 
Conversely, the Hamming metric can be defined in terms of the Hamming weight as 
$$
d_H(v,w) = \omega_H(vw^{-1}).
$$

Now, let us assume that $K$ is the additive group of the finite field $\mathbb{F}_q$. In \cite{BGL}, Brualdi, {
Graves, and Lawrence
} introduced metrics on the $\mathbb{F}_q$-vector space $K^N$ (which is isomorphic to $\mathbb{F}_q^N$) defined by posets. Let $P$ be a poset on the set $[N] = \{1, \dots, N\}$. The \textit{$P$-weight} of $v = (v_1, \dots, v_N) \in \mathbb{F}_q^N$ is defined by
$$
\omega_P(v) = |\{i \in [N] \mid i \preceq j \text{ for some } j \in \operatorname{supp}(v)\}|.
$$
Brualdi, {
Graves, and Lawrence
} showed that $d_P(v, w) = \omega_P(v - w)$ defines a metric on $\mathbb{F}_q^N$. 
This concept was generalized to digraphs in \cite{EFM2018}. 
If $\mathcal{G}$ is a digraph on $[N]$, the $\mathcal{G}$-\textit{weight} of a vector $v$ is defined as
$$
\omega_\mathcal{G}(v) = |\{i \in [N] \mid \text{there is a directed path, possibly trivial, from } i \text{ to some } j \in \operatorname{supp}(v)\}|.
$$

\subsection{Graphons.}

In precise terms, a {\em graphon} is a function 
$W: [0, 1] \times [0, 1] \to [0, 1]$ 
which satisfies the following two properties:
\begin{enumerate}
    \item[(1)] For all $ x, y \in [0, 1] $,
    \[
    W(x, y) = W(y, x).
    \]
    \item[(2)] The function $ W $ is measurable with respect to the Lebesgue measure on $[0, 1] \times [0, 1]$.
\end{enumerate}
The motivation behind this definition can be explained as follows. 
Suppose we have a finite graph $ \mathcal{G}(n) $ with $ n $ vertices, labeled as $ v_1, v_2, \dots, v_n $. 
The adjacency matrix $ A $ of $ \mathcal{G}(n) $ is an $ n \times n $ matrix where the $(i,j)$-th entry is given by 
\[
A_{ij} = 
\begin{cases} 
1 & \text{if there is an edge between } v_i \text{ and } v_j, \\
0 & \text{otherwise.}
\end{cases}
\]
To connect the adjacency matrix to a continuous function on the unit square $[0, 1] \times [0, 1]$, we divide $[0, 1]$ into $ n $ equal subintervals defined by 
\[
I_k := \left[\frac{k-1}{n}, \frac{k}{n}\right), \quad k = 1, 2, \dots, n.
\]
Each vertex $ v_k $ is mapped to the corresponding interval $ I_k $.

Over the unit square $[0, 1] \times [0, 1]$, we construct a piecewise constant function $ W_n(x, y) $ based on the adjacency matrix $ A $ by setting 
\[
W_n(x, y) := A_{ij} \quad \text{if } x \in I_i \text{ and } y \in I_j.
\]
In other words, the value of $ W_n(x, y) $ depends on which subintervals $ x $ and $ y $ belong to, and it is constant within each $ I_i \times I_j $ block.
As the number of vertices $ n $ increases, the subintervals $ I_k $ become smaller, and the step function $ W_n(x, y) $ becomes finer.
{
For a convergent sequence of dense graphs, the associated step functions converge to a graphon in cut distance, which allows relabeling of vertices. The limiting graphon is measurable and need not be continuous. Such convergence does not in general imply convergence in the $L^1$ norm.
}

For a comprehensive introduction to graphons, we recommend the book~\cite{Lovasz}.

\section{Proof of the main result}\label{S:Joint}

We begin with an elementary lemma.
\begin{Lemma}\label{L:subadditivefunction}
Let $f(t)$ denote the real-valued function defined by $f(t) = 1-e^{-t}$.
Then $f(t)$ is subadditive on the interval $[0,\infty)$.
In other words, we have 
$$
f(a + b) \leq f(a) + f(b)
$$ for 
every $a, b \geq 0${
.
}
\end{Lemma}

\begin{proof}
The inequality we want to prove is given by 
\begin{align*}
1 - e^{-(a+b)} &\leq (1 - e^{-a}) + (1 - e^{-b}).
\end{align*}
To get a clearer form of this inequality we rearrange the terms:
\begin{align*}
0 \leq 1 - e^{-a} - e^{-b} + e^{-a}e^{-b}.
\end{align*}
But the right hand side of this inequality is equal to $(1 - e^{-a})(1 - e^{-b})$.
Hence, the inequality we wish to prove is equivalent to the inequality 
$0 \leq (1 - e^{-a})(1 - e^{-b})$.
Now, since $a$ and $b$ are non-negative, $e^{-a}$ and $e^{-b}$ are both in the interval $(0, 1]$. 
Therefore, $(1 - e^{-a})$ and $(1 - e^{-b})$ are both in the interval $[0, 1]$. 
Thus, the inequality $(1 - e^{-a})(1 - e^{-b}) \geq 0$ always holds for $a, b \geq 0$.
This completes the proof of the fact that $f(t) = 1 - e^{-t}$ is a subadditive function on the entire interval $[0, \infty)$.
\end{proof}

The next two elementary lemmas will be useful for our next theorem. Transformations of metrics of this type are closely related to the classical work of Schoenberg~\cite{Schoenberg}, but the following direct proof is sufficient for our purposes.

\begin{Lemma}\label{L:power-subadditive}
Let $\alpha\geq 1$ and let $\phi_\alpha:[0,1]\to[0,1]$ be defined by
\[
\phi_\alpha(t)=1-(1-t)^\alpha.
\]
Then $\phi_\alpha$ is increasing and subadditive on $[0,1]$ in the following sense: if $u,v\in[0,1]$, then
\[
\phi_\alpha(\min\{1,u+v\})\leq \phi_\alpha(u)+\phi_\alpha(v).
\]
Consequently, if $d$ is a $[0,1]$-valued metric, then $\phi_\alpha\circ d$ is again a $[0,1]$-valued metric.
\end{Lemma}

\begin{proof}
The function $\phi_\alpha$ is increasing. If $u+v\leq 1$, then
fix $u$ and define
\[
H(v)=\phi_\alpha(u)+\phi_\alpha(v)-\phi_\alpha(u+v).
\]
Then $H(0)=0$ and
\[
H'(v)=\alpha\left((1-v)^{\alpha-1}-(1-u-v)^{\alpha-1}\right)\geq0
\]
for $0\leq v\leq1-u$, because $\alpha\geq1$ and $1-v\geq1-u-v$. Hence $H(v)\geq0$, proving
\[
\phi_\alpha(u+v)\leq \phi_\alpha(u)+\phi_\alpha(v).
\]
If $u+v\geq1$, then
\[
(1-u)+(1-v)\leq1.
\]
Since $\alpha\geq1$, we have $(1-u)^\alpha\leq 1-u$ and $(1-v)^\alpha\leq1-v$, and hence
\[
\phi_\alpha(u)+\phi_\alpha(v)=2-(1-u)^\alpha-(1-v)^\alpha\geq1=\phi_\alpha(1).
\]
The final assertion follows from the triangle inequality for $d$ together with monotonicity and subadditivity of $\phi_\alpha$.
\end{proof}

\begin{Lemma}\label{L:probabilistic-sum}
Let $\delta_1,\dots,\delta_m$ be $[0,1]$-valued metrics on sets $Y_1,\dots,Y_m$, respectively. Then
\[
\Delta((x_1,\dots,x_m),(y_1,\dots,y_m))
:=1-\prod_{r=1}^m\bigl(1-\delta_r(x_r,y_r)\bigr)
\]
is a $[0,1]$-valued metric on $Y_1\times\cdots\times Y_m$.
\end{Lemma}

\begin{proof}
It is enough to prove the assertion for two metrics, since the general case follows by induction. Put
$F(s,t)=s+t-st=1-(1-s)(1-t)$. If $a,b,c,d\in[0,1]$, then
\[
F(\min\{1,a+b\},\min\{1,c+d\})\leq F(a,c)+F(b,d).
\]
When $a+b\leq1$ and $c+d\leq1$, this follows from
\[
F(a,c)+F(b,d)-F(a+b,c+d)=ad+bc\geq0.
\]
If $a+b\geq1$ and $c+d\leq1$, then the left hand side is $1$, while
\[
F(a,c)+F(b,d)\geq F(a,c)+F(1-a,d)=1+c(1-a)+ad\geq1.
\]
The case $a+b\leq1$ and $c+d\geq1$ is symmetric. If $a+b\geq1$ and $c+d\geq1$, then the left hand side is again $1$, and the monotonicity of $F$ gives
\[
F(a,c)+F(b,d)\geq F(a,c)+F(1-a,1-c)=1+a+c-2ac\geq1.
\]
Now let $x,y,z\in Y_1\times Y_2$. Applying this inequality with
$a=\delta_1(x_1,y_1)$, $b=\delta_1(y_1,z_1)$, $c=\delta_2(x_2,y_2)$, and $d=\delta_2(y_2,z_2)$ gives the triangle inequality for $F(\delta_1,\delta_2)$.
The remaining metric axioms are immediate.
\end{proof}

Let $(X_1, d_1),\dots, (X_N, d_N)$ be metric spaces, where {each $d_i$ takes values in $[0,1)$}. 
For positive real numbers $a_1,\dots, a_N$, we define a function $d_{a_1,\dots,a_N}$ on the product $X_1 \times \cdots \times X_N$ by setting 
\begin{align}\label{A:das}
d_{a_1,\dots,a_N}(\mathbf{x},\mathbf{y}) := 1 - \prod_{i=1}^N (1 - d_i(x_i, y_i))^{a_i},
\end{align}
where $\mathbf{x}=(x_1,\dots, x_N)$ and $\mathbf{y}=(y_1,\dots, y_N)$ are from $X_1\times \cdots \times X_N$.

\begin{Theorem}\label{T:firstmaintheorem}
Let $(a_1,\dots, a_N)$ be a list of real numbers such that $a_i\geq 1$ for $i=1,\dots, N$. 
Let {
$d_{a_1,\dots, a_N} : (X_1\times \cdots \times X_N)^2\to \R$
} be the function defined above. 
Then $d_{a_1,\dots, a_N}$ is a metric that takes values in $[0,1)$.
\end{Theorem}

\begin{proof}
For each $i$, Lemma~\ref{L:power-subadditive} shows that
\[
\delta_i(x_i,y_i):=1-(1-d_i(x_i,y_i))^{a_i}
\]
is a $[0,1]$-valued metric on $X_i$. Since
\[
1-d_{a_1,\dots,a_N}(\mathbf{x},\mathbf{y})
=\prod_{i=1}^N\bigl(1-\delta_i(x_i,y_i)\bigr),
\]
Lemma~\ref{L:probabilistic-sum} implies that $d_{a_1,\dots,a_N}$ is a metric on the product $X_1\times\cdots\times X_N$. Its values lie in $[0,1)$ because the component distances take values in $[0,1)$.
\end{proof}

We will strengthen Theorem~\ref{T:firstmaintheorem} by relaxing boundaries.
More precisely, we will allow the distance functions $d_1,\dots, d_N$ {
to take values
} in $[0, 1]$ instead of $[0,1)$.

\begin{Theorem}\label{T:firstmaintheorem.1}
Let $(a_1,\dots, a_N)$ be a list of real numbers such that $a_i\geq 1$ for $i=1,\dots, N$. 
Let {
$d_{a_1,\dots, a_N} : (X_1\times \cdots \times X_N)^2\to \R$
} be the function defined in (\ref{A:das}).
Then $d_{a_1,\dots, a_N}$ is a metric that takes values in $[0,1]$.
\end{Theorem}

\begin{proof}
The fact that $d_{a_1,\dots, a_N}$ takes values in $[0,1]$ is evident. 
The non-negativity, identity of indiscernibles, and symmetry properties of the function $d_{a_1,\dots,a_N}$ are straightforward to check. 
It remains to prove that the function $d_{a_1,\dots,a_N}$ satisfies the triangle inequality.

{
There is no singularity in the product formula when a coordinate distance equals $1$, since its factor is $0^{a_i}=0$ and $a_i\geq1$.
If this happens for $(\mathbf x,\mathbf y)$, then the product distance is $1$.
The coordinate triangle inequality forces $d_i(x_i,z_i)+d_i(z_i,y_i)\geq1$.
Also, $1-(1-t)^{a_i}\geq t$ on $[0,1]$, and the full product distance is at least this transformed coordinate distance.
Thus the two distances on the right side of the triangle inequality have sum at least $1$.
If a distance on the right already equals $1$, the inequality follows from the bound $d_{a_1,\dots,a_N}\leq1$.
The following cases give the details.
}

Let $\mathbf{x}=(x_1,\dots, x_N)$, $\mathbf{y}=(y_1,\dots, y_N)$, and $\mathbf{z}=(z_1,\dots, z_N)$ be three $N$-tuples from $X_1\times \cdots \times X_N$.  

If none of the distances $d_i(x_i,y_i),d_i(y_i,z_i)$, and $d_i(x_i,z_i)$, where $i=1,\dots, N$, is 1, 
then the proof of Theorem~\ref{T:firstmaintheorem} shows that $d_{a_1,\dots, a_N}(\mathbf{x},\mathbf{y})\leq d_{a_1,\dots, a_N}(\mathbf{x},\mathbf{z})+d_{a_1,\dots, a_N}(\mathbf{z},\mathbf{y})$.
We proceed with the assumption that, for some $i\in \{1,\dots, N\}$, one of the numbers $d_i(x_i,y_i),d_i(y_i,z_i)$, or $d_i(x_i,z_i)$ is 1.

\textbf{Case 1.}
We assume that $d_i(x_i,y_i)=1$ for some $i\in \{1,\dots, N\}$. 
We set $i=1$ for simplicity.
Under this assumption, we will use mathematical induction on $N$ to prove that $d_{a_1,\dots, a_N}$ satisfies the triangle inequality. 
We begin with the base case, where $N=1$.
Then our goal is to show the following inequality
\begin{align*}
 1- (1-a)^{a_1} \leq \left( 1- (1-b)^{a_1}\right) + \left(1- (1-c)^{a_1}\right),
\end{align*}
where $a=d_1(x_1,y_1)$, $b=d_1(x_1,z_1)$, and $c=d_1(z_1,y_1)$ are such that $a=1$, $1\leq b+c$, and $0\leq b,c\leq 1$. 
Equivalently, we will show that 
\begin{align}\label{A:beforesimplifying}
(1-b)^{a_1} +  (1-c)^{a_1} \leq 1, 
\end{align}
where we have $1\leq b+c$ and $0\leq b,c\leq 1$. 
But these constraints describe the upper right triangle, denoted by $B$, in the unit square in Figure~\ref{F:shaded}.
\begin{figure}[htp]
\centering
\scalebox{0.75}{  
\begin{tikzpicture}
    \draw[thick] (0, 0) rectangle (4, 4);

    \draw[thick, red] (0, 4) -- (4, 0);

    \fill[blue!20] (0, 4) -- (4, 0) -- (4, 4) -- cycle;

    \node at (-0.3, -0.3) {$0$};
    \node at (4.3, -0.3) {$1$};
    \node at (-0.3, 4.3) {$1$};
    \node at (4.3,4.3) {$(1,1)$};

    \draw[->] (-0.5, 0) -- (4.5, 0) node[right] {$b$};
    \draw[->] (0, -0.5) -- (0, 4.5) node[above] {$c$};

    \node[red] at (3, 1) {$b + c = 1$};

\node at (1,1) {$A$};
\node at (3,3) {$B$};

\end{tikzpicture}}
\caption{The shaded triangle $B$ is where $0\leq b,c\leq 1$ and $b+c\geq 1$ hold.}
\label{F:shaded}
\end{figure}
Since we are interested in the inequality (\ref{A:beforesimplifying}), we apply the transformation $(b,c)\mapsto (1-b,1-c)$
to triangle $B$.
The image of $B$ under this transformation is precisely the non-shaded triangle, denoted by $A$, in Figure~\ref{F:shaded}.
On the triangle $A$ we have $x,y\geq0$ and $x+y\leq1$. Since $\alpha\geq1$ implies $x^\alpha\leq x$ and $y^\alpha\leq y$ for $x,y\in[0,1]$, we have $x^\alpha+y^\alpha\leq x+y\leq1$. Taking $x=1-b$, $y=1-c$, and $\alpha=a_1$ proves (\ref{A:beforesimplifying}). This completes the proof of the base case of our induction.

We now assume that our claim holds for $N-1$, and proceed to prove it for $N$. 
Notice the decomposition 
\begin{align}\label{A:das1}
d_{a_1,\dots,a_N}(\mathbf{x},\mathbf{y}) = {\underbrace{1-(1-d_1(x_1,y_1))^{a_1}}_{\alpha(\mathbf{x},\mathbf{y})}}+
{\underbrace{(1-d_1(x_1,y_1))^{a_1}}_{\beta(\mathbf{x},\mathbf{y})}} \underbrace{ \left(1- \prod_{i=2}^N (1 - d_i(x_i, y_i))^{a_i}\right)}_{\gamma(\mathbf{x},\mathbf{y})},
\end{align}
where $\alpha(\mathbf{x},\mathbf{y})$, $\beta(\mathbf{x},\mathbf{y})$, and $\gamma(\mathbf{x},\mathbf{y})$ are nonnegative numbers.
In fact, under our assumptions, we have $\alpha(\mathbf{x},\mathbf{y})=1$, $\beta(\mathbf{x},\mathbf{y})=0$, and $\gamma(\mathbf{x},\mathbf{y})\geq 0$. 
We have similar decompositions for $d_{a_1,\dots,a_N}(\mathbf{x},\mathbf{z})$ and $d_{a_1,\dots,a_N}(\mathbf{z},\mathbf{y})$.
Our goal is to prove the triangle inequality, 
\begin{align}\label{A:triangleineq}
1 \leq \left(\alpha(\mathbf{x},\mathbf{z}) + \beta(\mathbf{x},\mathbf{z})\gamma(\mathbf{x},\mathbf{z}) \right) +\left(\alpha(\mathbf{z},\mathbf{y}) + \beta(\mathbf{z},\mathbf{y})\gamma(\mathbf{z},\mathbf{y}) \right). 
\end{align}
But we already proved in the base case earlier that the inequality $1\leq \alpha(\mathbf{x},\mathbf{z}) + \alpha(\mathbf{z},\mathbf{y})$ holds. 
Since $\beta(\mathbf{x},\mathbf{z})\gamma(\mathbf{x},\mathbf{z}) +  \beta(\mathbf{z},\mathbf{y})\gamma(\mathbf{z},\mathbf{y})$ is nonnegative, the inequality in (\ref{A:triangleineq}) follows. 
This completes the proof of our assertion that if $d_1(x_1,y_1)=1$, then $d_{a_1,\dots,a_N}(\mathbf{x},\mathbf{y})\leq d_{a_1,\dots,a_N}(\mathbf{x},\mathbf{z})+d_{a_1,\dots,a_N}(\mathbf{z},\mathbf{y})$ holds. 
\medskip

\textbf{Case 2.}
We assume that one of the numbers $d_i(x_i,z_i)$ or $d_i(z_i,y_i)$ is 1 for some $i\in \{1,\dots, N\}$. 
Since there is no difference between our arguments for $d_i(x_i,z_i)=1$ or $d_i(z_i,y_i)=1$, we proceed with the assumption
that $d_1(x_1,z_1)=1$. 
Hence, we know that $d_{a_1,\dots, a_N}(\mathbf{x},\mathbf{z}) =1$. 
But we already pointed out at the beginning of our proof that $d_{a_1,\dots, a_N}$ takes values in $[0,1]$.
Therefore, the inequality $d_{a_1,\dots, a_N}(\mathbf{x},\mathbf{y}) \leq 1+ d_{a_1,\dots, a_N}(\mathbf{z},\mathbf{y})$ holds. 
This completes the proof of our theorem.
\end{proof}

We are now ready to prove the main result, which we restate here for the sake of convenience.

\medskip

Let $(X_1,d_1),\dots,(X_N,d_N)$ be metric spaces whose distance functions $d_i$, $i=1,\dots,N$, take values in the interval $[0,1]$.
Let $\mathcal{X}$ denote the set theoretic product $X_1\times \cdots \times X_N$.
Let $\mathcal{G}$ be a directed graph with vertex set $\{1,\dots,N\}$.
For each $j\in[N]$, define the closed in-neighborhood
\[
B_j(\mathcal{G})=\{j\}\cup\{i:(i,j)\in E(\mathcal{G})\}.
\]
If every edge $(i,j)$ of $\mathcal{G}$ is assigned a probability score $p_{ij}\in(0,1]$, and if we set $p_{jj}=1$ for every $j$, then the mapping
\begin{align}\label{A:initialformofd}
d_{\mathcal{X},\mathcal{G},\mathcal{P}}(\mathbf{g},\mathbf{h})
:=1-\frac{1}{N}\sum_{j=1}^N \prod_{i\in B_j(\mathcal{G})}\left[1-d_i(g_i,h_i)\right]^{1/p_{ij}}
\end{align}
for $\mathbf{g}=(g_1,\dots,g_N)$ and $\mathbf{h}=(h_1,\dots,h_N)$ defines a normalized metric on $\mathcal{X}$.

\begin{proof}[Proof of the Main Theorem]
For $j=1,\dots,N$, define
\[
d_{P,j}(\mathbf{x},\mathbf{y})
:=1-\prod_{i\in B_j(\mathcal{G})}\left(1-d_i(x_i,y_i)\right)^{1/p_{ij}},
\]
where $\mathbf{x}=(x_1,\dots,x_N)$ and $\mathbf{y}=(y_1,\dots,y_N)$ are from $\mathcal{X}$. Since $p_{ij}\in(0,1]$, we have $1/p_{ij}\geq1$.
By Theorem~\ref{T:firstmaintheorem.1}, the function $d_{P,j}$ is a metric on the subproduct indexed by $B_j(\mathcal{G})$ and hence a pseudometric on the full product $\mathcal{X}$.
Moreover,
\[
d_{\mathcal{X},\mathcal{G},\mathcal{P}}(\mathbf{x},\mathbf{y})=\frac1N\sum_{j=1}^N d_{P,j}(\mathbf{x},\mathbf{y}).
\]
A nonnegative linear combination of pseudometrics is again a pseudometric; hence $d_{\mathcal{X},\mathcal{G},\mathcal{P}}$ is nonnegative, symmetric, and satisfies the triangle inequality.
It remains to verify identity of indiscernibles. If $\mathbf{x}\neq\mathbf{y}$, choose $k$ such that $x_k\neq y_k$. Since $k\in B_k(\mathcal{G})$ and $p_{kk}=1$, the $k$-th summand satisfies
\[
d_{P,k}(\mathbf{x},\mathbf{y})
=1-\prod_{i\in B_k(\mathcal{G})}\left(1-d_i(x_i,y_i)\right)^{1/p_{ik}}
\geq 1-(1-d_k(x_k,y_k))>0.
\]
Thus $d_{\mathcal{X},\mathcal{G},\mathcal{P}}(\mathbf{x},\mathbf{y})>0$. This completes the proof.
\end{proof}

\section{Proof of Corollary~\ref{intro:C1} and a sphere packing bound}\label{S:Specializations}

In this section we use the unnormalized distance
\[
D_{X,G,P}:=N d_{X,G,P}.
\]
This is the natural scale for coding-theoretic applications, since in the binary specialization it is integer-valued.

We first prove Corollary~\ref{intro:C1}. Let us recall its statement for convenience. 
\medskip

We assume that the underlying sets in the main result, $X_1,\dots, X_N$, are finite and that the corresponding metrics $d_1,\dots, d_N$ take values in $\{0,1\}$. Let $G^{\preceq}$ denote the oriented reflexive transitive closure of a DAG $G$, oriented so that the closed in-neighborhood of $j$ is its forward reachability set:
\[
B_j(G^{\preceq})=\{i\in[N]: \text{there is a directed path, possibly trivial, from } j \text{ to } i \text{ in }G\}.
\]
Then, for every $\mathbf{g}=(g_1,\dots,g_N) \in X$ and $\mathbf{h} =(h_1,\dots,h_N)\in X$, we have 
\begin{align*}
D_{X,G^{\preceq},P}(\mathbf{g},\mathbf{h})=|\langle {\rm supp}(\mathbf{g},\mathbf{h})\rangle|,
\qquad
 d_{X,G^{\preceq},P}(\mathbf{g},\mathbf{h})=\frac{|\langle {\rm supp}(\mathbf{g},\mathbf{h})\rangle|}{N},
\end{align*}
where ${\rm supp} (\mathbf{g},\mathbf{h})=\{j\in[N]:g_j\neq h_j\}$ and $\langle {\rm supp}(\mathbf{g},\mathbf{h})\rangle$ is the hereditary closure of the support.

\begin{proof}[Proof of Corollary~\ref{intro:C1}] 
From the definition of the unnormalized distance,
\[
D_{X,G^{\preceq},P}(\mathbf{g},\mathbf{h})
=\sum_{j=1}^N\left(1-\prod_{i\in B_j(G^{\preceq})}\left(1-d_i(g_i,h_i)\right)^{1/p_{ij}}\right).
\]
Since each $d_i(g_i,h_i)$ is either $0$ or $1$, the $j$-th summand is equal to $1$ exactly when
\[
B_j(G^{\preceq})\cap {\rm supp}(\mathbf{g},\mathbf{h})\neq\emptyset,
\]
and it is equal to $0$ otherwise. By the definition of $G^{\preceq}$, this condition says precisely that there is a directed path in $G$ from $j$ to an element of ${\rm supp}(\mathbf{g},\mathbf{h})$. Equivalently, $j$ belongs to the hereditary closure $\langle {\rm supp}(\mathbf{g},\mathbf{h})\rangle$. Therefore the sum counts the elements of this hereditary closure. Dividing by $N$ gives the normalized formula for $d_{X,G^{\preceq},P}$.
\end{proof}

We are now ready to prove a sphere-packing bound for our metrics. 
Let us recall its statement. 

\medskip

Let $X=\F_q^N$, $d_i\in\{0,1\}$, and suppose $G$ is a DAG. For $S\subseteq[N]$, write $\langle S\rangle$ for the smallest lower (hereditary) set containing $S$ in the poset induced by $G$, and define the $G$-weight $\wtG(x):=|\langle \mathrm{supp}(x)\rangle|$. For $r\in\{0,1,\dots,N\}$ set the (closed) ball
\begin{equation*}
\ballG(r)\;:=\;\{x\in X:\ \wtG(x)\le r\}.
\end{equation*}
If a code $C\subset X$ has minimum unnormalized $G$-distance $D_{\min}$ with respect to $D_{X,G^{\preceq},P}$, then
\begin{equation}
|C|
\;\le\;\frac{|X|}{\displaystyle\max_{z\in X}\bigl|\,z+\ballG(\lfloor (D_{\min}-1)/2\rfloor)\,\bigr|}
\;\le\;\frac{q^N}{\displaystyle\sum_{s=0}^{\lfloor (D_{\min}-1)/2\rfloor} J_G(s)\,(q-1)^s}\,,
\end{equation}
where $J_G(s)$ denotes the number of lower ideals of cardinality $s$ in the poset of $G$.

\begin{proof}[Proof of Theorem~\ref{thm:poset-sphere-packing}]
Let $t:=\lfloor(D_{\min}-1)/2\rfloor$. Since $D_{X,G^{\preceq},P}(x,y)=\wtG(x-y)$ is a translation-invariant metric,
the closed balls $\{z+\ballG(t):z\in C\}$ are pairwise disjoint. Indeed, if the intersection 
$(z_1+\ballG(t))\cap(z_2+\ballG(t))$ is nonempty, then for a point $x$ in the intersection we find 
$$
D_{X,G^{\preceq},P}(z_1,z_2)\le D_{X,G^{\preceq},P}(z_1,x)+D_{X,G^{\preceq},P}(x,z_2)\le t+t=2t<D_{\min},
$$ 
which is absurd. 
Therefore, the first claimed inequality readily holds: 
\[
|C|\cdot|\ballG(t)| \;\le\; |X|=q^N
\qquad\Longrightarrow\qquad
|C|\;\le\;\frac{q^N}{|\ballG(t)|}.
\]

It remains to give the lower bound on $|\ballG(t)|$ used in \eqref{eq:packing-bound}.
For any lower order ideal $I\subseteq[N]$ of size $s\le t$, every vector $x$ with
$\mathrm{supp}(x)=I$ satisfies 
$$
\wtG(x)=|\langle\mathrm{supp}(x)\rangle|=|I|=s\le t.
$$ 
Hence, $x\in\ballG(t)$.
The number of vectors with $\mathrm{supp}(x)=I$ is exactly $(q-1)^{|I|}$. Summing over all lower ideals
of size $s\le t$ gives
\[
|\ballG(t)| \;\ge\; \sum_{s=0}^{t} J_G(s)\,(q-1)^s,
\]
where $J_G(s)$ counts lower ideals of size $s$ in the poset of $G$ (with $J_G(0)=1$ accounting for $x=0$).
Combining the two equations yields
\[
|C|\;\le\;\frac{q^N}{|\ballG(t)|}\;\le\;\frac{q^N}{\sum_{s=0}^{t} J_G(s)\,(q-1)^s},
\]
as claimed.
\end{proof}

\section{Proof of Corollary~\ref{intro:C2}}\label{S:Graphons}

Recall the statement of Corollary~\ref{intro:C2}. { We retain its assumption that $\rho(u,v):=-\log(1-d(u,v))$ is a metric.} In the limit $N \rightarrow \infty$, where the graph $G$ in the main result approaches a graphon ${W:[0,1]^2 \rightarrow (0,1]}$, defined over the random variables $x,y \in [0,1]$, then
{
\begin{align*}
d_{G,W} (g,h) := 1 - \mathbb{E}_{x}\left[ \exp \left(\mathbb{E}_{y} \left[\frac{\log \left(1- d(g(y),h(y))\right)}{W(x,y)} \right]\right) \right]
\end{align*}
is a metric defined over the {continuous} functions $g,h:[0,1] \rightarrow \mathbb{C}$, where
$\mathbb{E}_{x}(\cdot)$ and $\mathbb{E}_{y}(\cdot)$ denote expectations with respect to $x$ and $y$, respectively, whenever the expectations are well defined.
}
\medskip

{
\begin{proof}[Proof of Corollary~\ref{intro:C2}]
We first modify the original definition of $d_{X,G,P} (g,h)$ as follows: 
\begin{align}\label{A:initialformofdcond}
d_{X,G,P} (g,h) := 1 - \frac{1}{N}\sum_{j=1}^N \exp \left[\sum_{i=1}^N \frac{\log \left(1- d(g_j,h_j)\right)}{H_{ji}p_{j|i}}\right].
\end{align}
where the indices $i,j$ represent discrete random variables drawn from a joint probability mass function $p_{ji} \ge 0,\sum_{i,j} p_{j,i} = 1$ and $p_{j|i} \ge 0$ represents a conditional probability measure satisfying $\sum_j^N p_{j|i} = 1$. $H_{ji} \in (0,1]$ in equation~\ref{A:initialformofdcond} represents an auxiliary variable indexed by the random variables $i$ and $j$. { Note that for this corollary, we restrict $d(g_j,h_j) < 1, \forall j$ which ensures that the logarithm and the exponential functions in equation~\ref{A:initialformofdcond} are well defined.} {We rewrite} equation~\ref{A:initialformofdcond} in terms of the joint and the marginal probabilities as
\begin{align*}
d_{X,G,P} (g,h) := 1 - \frac{1}{N}\sum_{j=1}^N \exp \left[\sum_{i=1}^N \frac{p_i}{H_{ji}p_{ji}}\log \left(1- d(g_j,h_j)\right)\right].
\end{align*}
Denoting a discrete graphon function $W_{ji} = H_{ji}p_{ji}$, the preceding display leads to
\begin{align*}
d_{X,G,P} (g,h) := 1 - \frac{1}{N}\sum_{j=1}^N \exp \left[\sum_{i=1}^N \frac{p_i}{W_{ji}}\log \left(1- d(g_j,h_j)\right)\right]
\end{align*}
where the metric can be viewed as a distance function that is estimated using finite $N \times N$ samples drawn from an underlying distribution. In the limit $N \rightarrow \infty$, and similar to the framework of graphons, where the vertices $i,j$ can be associated with random variables $x,y \in \left[0,1\right]$ with marginal distributions $p(x),p(y)$. 
Therefore, as $N \rightarrow \infty$, the distance $d_{X,G,P} (g,h) \rightarrow d_{X,G,W}$ can be extended to continuous functions $g,h: \left[0,1\right] \rightarrow \mathbb{C}$ over graphons as 
\begin{align*}
d_{G,W} (g,h)  := 1 - \int_0^1 \exp \left[\int_0^1 \frac{\log \left(1- d(g(y),h(y))\right)}{W(x,y)}p(y) dy\right] p(x) dx.
\end{align*}
{ Since the $\exp(.)$ function in equation~\ref{A:initialformofdcond} is bounded, for the last step we have used: (a) convergence in the continuous limit for the {Lebesgue} integration; and (b) the law of large numbers to replace empirical mean by the true expectation.} This finishes the proof of our corollary. 
\end{proof}
}

{
For a simple example, take $W(x,y)=w\in(0,1]$, let $x$ and $y$ be independent and uniform on $[0,1]$, and set $d(u,v)=1-e^{-|u-v|}$.
The formula becomes
\[
d_{G,W}(g,h)=1-\exp\left(-\frac1w\int_0^1|g(y)-h(y)|\,dy\right).
\]
For $g(y)=0$ and $h(y)=y$, this gives $1-e^{-1/(2w)}$.
In numerical evaluation, small values of $W(x,y)$ or values of $d(g(y),h(y))$ close to $1$ require care because the exponent contains $1/W(x,y)$ and $\log(1-d(g(y),h(y)))$.
}

\section{Proofs of Corollaries~\ref{intro:C3}, \ref{intro:C4}, and Proposition~\ref{P:cartesian-factorization}}\label{S:Semiring}

{
Semirings generalize rings
} by dropping the requirement for additive inverses. This allows for algebraic structures in various fields where subtraction does not naturally exist or is not meaningful, such as in the context of optimization, theoretical computer science, and indeed, graph theory.
For a broad review of the theory of semirings and their applications to these fields, we recommend the monograph~\cite{GondranMinoux}.

In this section, we discuss the properties of joint metrics with respect to semiring structure on directed graphs. We begin with the proof of Corollary~\ref{intro:C3}. 
We recall its statement for convenience. 
\bigskip

Let $(X_1,d_1),\dots, (X_N,d_N)$ be metric spaces as in the main result of our paper.
Let $\mathcal{G}_1$ (resp. $\mathcal{G}_2$) be a directed graph on the set $\{1,\dots, N_1\}$ (resp. on the set $\{N_1+1,\dots, N\}$). 
We denote by $\mathcal{X}^{(1)}$ (resp. by $\mathcal{X}^{(2)}$) the product metric space $X_1\times \cdots \times X_{N_1}$ (resp. $X_{N_1+1}\times \cdots \times X_N$).
Let $\mathcal{P}_1$ (resp. $\mathcal{P}_2$) be the matrix of probability scores on the edges of $\mathcal{G}_1$ (resp. on the edges of $\mathcal{G}_2$). 
If the directed graph $\mathcal{G}$ is a disjoint union of two sub-directed graphs, $\mathcal{G} := \mathcal{G}_1 \sqcup \mathcal{G}_2$,
then we have 
$$
d_{\mathcal{X},\mathcal{G},\mathcal{P}}\ = \ \frac{N_1}{N}\ d_{\mathcal{X}^{(1)}, \mathcal{G}_1 ,\mathcal{P}_1} \ + \ \frac{N_2}{N}\ d_{\mathcal{X}^{(2)}, \mathcal{G}_2 ,\mathcal{P}_2}.
$$
where $\mathcal{X} = \mathcal{X}^{(1)}\times \mathcal{X}^{(2)}$ and $N_2 = N-N_1$.

\begin{proof}[Proof of Corollary~\ref{intro:C3}]
Let $E_1$ (resp. $E_2$) denote the edge set of $\mathcal{G}_1$ (resp. of $\mathcal{G}_2$).  
We now look closely at the definition of $d_{\mathcal{X}, \mathcal{G}_1 \sqcup \mathcal{G}_2, \mathcal{P}}$.
Let $(\mathbf{g}, \mathbf{h})\in \mathcal{X}\times \mathcal{X}$. 
If $\mathbf{g}$ is given by $(g_1,\dots, g_N)$, then we will write $\mathbf{g}^{(1)}$ (resp. $\mathbf{g}^{(2)}$) for $(g_1,\dots, g_{N_1})$
(resp. for $(g_{N_1+1},\dots, g_N)$). 
The vectors $\mathbf{h}^{(1)}$ and $\mathbf{h}^{(2)}$ are defined similarly. 
Then we have 
\begin{align*}
d_{\mathcal{X}, \mathcal{G}_1 \sqcup \mathcal{G}_2, \mathcal{P}}(\mathbf{g}, \mathbf{h})
&= 1 - \frac{1}{N} \sum_{j=1}^N \prod_{i\in B_j(\mathcal{G}_1\sqcup\mathcal{G}_2)} \left[1 - d_i(g_i, h_i)\right]^{1/p_{ij}} \\
&= {1 - \frac{1}{N} \Biggl[\underbrace{ \left( \sum_{j=1}^{N_1} \prod_{i\in B_j(\mathcal{G}_1)} \left[1 - d_i(g_i, h_i)\right]^{1/p_{ij}} \right)}_{A}}\\
&\hspace{4em}{+\underbrace{\left( \sum_{j=N_1+1}^{N} \prod_{i\in B_j(\mathcal{G}_2)} \left[1 - d_i(g_i, h_i)\right]^{1/p_{ij}} \right)}_{B}\Biggr].}
\end{align*}
We reorganize this expression as follows: 
\begin{align*}
d_{\mathcal{X}, \mathcal{G}_1 \sqcup \mathcal{G}_2, \mathcal{P}}(\mathbf{g}, \mathbf{h}) 
&= 1 - \frac{1}{N} \left[ \left( -N_1 (-1 + 1 - \frac{1}{N_1} A) \right) +  \left( -N_2 (-1 + 1 - \frac{1}{N_2} B)  \right) \right] \\
&= 1 - \frac{1}{N} \left[ \left( N_1 - N_1\left(1 - \frac{1}{N_1} A\right) \right) +  \left( N_2-N_2 \left(1 - \frac{1}{N_2} B\right)  \right) \right] \\
&= 1-  \frac{N_1 + N_2}{N}  - \frac{1}{N} \left[ \left( - N_1\left(1 - \frac{1}{N_1} A\right) \right) +  \left(-N_2 \left(1 - \frac{1}{N_2} B\right)  \right) \right] \\
&=  \frac{N_1}{N} \left(1 - \frac{1}{N_1} A\right) +   \frac{N_2}{N} \left(1 - \frac{1}{N_2} B\right).
\end{align*}
Since the term $1 - \frac{1}{N_1} A$ (resp. the term $1 - \frac{1}{N_2} B$) is the metric value $d_{\mathcal{X}^{(1)}, \mathcal{G}_1 ,\mathcal{P}_1}(\mathbf{g}^{(1)},\mathbf{h}^{(1)})$ 
(resp. $d_{\mathcal{X}^{(2)}, \mathcal{G}_2 ,\mathcal{P}_2}(\mathbf{g}^{(2)},\mathbf{h}^{(2)})$),
the proof follows.
\end{proof}

Let $(X_1, d_1), \dots, (X_N, d_N)$ be metric spaces as in the main result of our paper. 
Consider a list of directed graphs $\mathcal{G}_1, \dots, \mathcal{G}_r$, where $\mathcal{G}_i$ has $N_i$ nodes for $i = 1, \dots, r$, and assume that $N_1 + \cdots + N_r = N$.  
Define the following subsets of the product space:  
\begin{itemize}  
    \item Let $\mathcal{X}^{(1)}$ denote the product of the first $N_1$ metric spaces, that is,  
    \[  
    \mathcal{X}^{(1)} = X_1 \times X_2 \times \cdots \times X_{N_1}.  
    \]  
    \item Similarly, let $\mathcal{X}^{(2)}$ denote the product of the next $N_2$ metric spaces, that is,  
    \[  
    \mathcal{X}^{(2)} = X_{N_1+1} \times X_{N_1+2} \times \cdots \times X_{N_1+N_2}.  
    \]  
    \item Continuing inductively, define $\mathcal{X}^{(i)}$ for $i = 3, \dots, r$ as the product of the next $N_i$ metric spaces, ensuring the main product $ X_1\times \cdots \times X_N$ is partitioned into the subproducts $\mathcal{X}^{(1)}, \mathcal{X}^{(2)}, \dots, \mathcal{X}^{(r)}$.  
\end{itemize}  
This construction associates each graph $\mathcal{G}_i$ with the corresponding product $\mathcal{X}^{(i)}$.  

In this notation, we have the following corollary of our previous result. 

\begin{Corollary}\label{C:introT11}
We maintain the notation from the previous paragraph. 
Let $\mathcal{X}$ denote the metric space $X_1\times \cdots \times X_N$.
Let $\mathcal{G}$ denote the directed graph obtained by taking the disjoint union of the graphs $\mathcal{G}_1,\dots, \mathcal{G}_r$. 
Finally, let $\mathcal{P}$ denote the block diagonal matrix obtained by putting together the matrices of probability scores $\mathcal{P}_1,\dots, \mathcal{P}_r$ of the edges of the graphs $\mathcal{G}_1,\dots, \mathcal{G}_r$ 
in the written order. 
Then we have 
$$
d_{\mathcal{X}, \mathcal{G} ,\mathcal{P}}\ = \ \frac{N_1}{N}\ d_{\mathcal{X}^{(1)}, \mathcal{G}_1 ,\mathcal{P}_1} \ +\cdots + \ \frac{N_r}{N}\ d_{\mathcal{X}^{(r)}, \mathcal{G}_r ,\mathcal{P}_r}.
$$
\end{Corollary}

\begin{proof}
In light of our previous theorem, the proof follows by induction on the number $r$ of directed graphs.
We omit details.  
\end{proof}

\bigskip

We now proceed to discuss what happens to our metrics under Cartesian products of directed graphs.
Let us recall the statement of our relevant result, Corollary~\ref{intro:C4}. 
\medskip

We maintain our notation from the introductory section. The exact formula stated in Corollary~\ref{intro:C4} expresses $d_{F,\mathcal{G}_1\square\mathcal{G}_2,\mathcal{P}}$ as an average of products over the vertex set of $\mathcal{G}_1\square\mathcal{G}_2$; unlike the disjoint-union case, this expression does not factor as a product of the two graph metrics attached to $\mathcal{G}_1$ and $\mathcal{G}_2$.
\medskip

\begin{proof}[Proof of Corollary~\ref{intro:C4}]
Let $(\mathbf{g},\mathbf{h})\in F\times F$ and write
\[
\delta_{u_1,u_2}:=d_{u_1,N_1+u_2}(g_{u_1,u_2},h_{u_1,u_2}).
\]
For a vertex $(v_1,v_2)$ of $\mathcal{G}_1\square\mathcal{G}_2$, its closed in-neighborhood is
\[
\{(v_1,v_2)\}
\cup \{(u_1,v_2):(u_1,v_1)\in E_1\}
\cup \{(v_1,u_2):(u_2,v_2)\in E_2\}.
\]
Substituting this closed in-neighborhood into the definition of $d_{F,\mathcal{G}_1\square\mathcal{G}_2,\mathcal{P}}$ gives
\begin{align*}
d_{F, \mathcal{G}_1 \square \mathcal{G}_2, \mathcal{P}}(\mathbf{g}, \mathbf{h})
=1-\frac{1}{N_1N_2}\sum_{v_1\in V_1}\sum_{v_2\in V_2}
&\left(1-\delta_{v_1,v_2}\right)\\
&\times\prod_{(u_1,v_1)\in E_1}\left(1-\delta_{u_1,v_2}\right)^{1/p^{(1)}_{u_1v_1}}\\
&\times\prod_{(u_2,v_2)\in E_2}\left(1-\delta_{v_1,u_2}\right)^{1/p^{(2)}_{u_2v_2}},
\end{align*}
which is the claimed formula. The last assertion follows because, for general matrices of distances $\delta_{u_1,u_2}$, an average of products cannot be rewritten as the product of the corresponding row and column averages without additional separability assumptions.
\end{proof} 

\begin{proof}[Proof of Proposition~\ref{P:cartesian-factorization}]
Starting from the exact formula in Corollary~\ref{intro:C4}, write
\[
1-d_{F,\mathcal{G}_1\square\mathcal{G}_2,\mathcal{P}}(\mathbf{g},\mathbf{h})
=\frac{1}{N_1N_2}\sum_{v_1\in V_1}\sum_{v_2\in V_2}
s_{v_1,v_2}
\prod_{(u_1,v_1)\in E_1}s_{u_1,v_2}^{1/p^{(1)}_{u_1v_1}}
\prod_{(u_2,v_2)\in E_2}s_{v_1,u_2}^{1/p^{(2)}_{u_2v_2}}.
\]
Substituting $s_{u_1,u_2}=r_{u_1}c_{u_2}$, the summand indexed by $(v_1,v_2)$ becomes
\[
r_{v_1}c_{v_2}
\prod_{(u_1,v_1)\in E_1}
(r_{u_1}c_{v_2})^{1/p^{(1)}_{u_1v_1}}
\prod_{(u_2,v_2)\in E_2}
(r_{v_1}c_{u_2})^{1/p^{(2)}_{u_2v_2}}.
\]
Using $A(v_1)=A$ and $B(v_2)=B$, this factors as
\[
\left[
r_{v_1}^{1+B}
\prod_{(u_1,v_1)\in E_1}
r_{u_1}^{1/p^{(1)}_{u_1v_1}}
\right]
\left[
c_{v_2}^{1+A}
\prod_{(u_2,v_2)\in E_2}
c_{u_2}^{1/p^{(2)}_{u_2v_2}}
\right].
\]
Therefore the double average factors as $S_1(r)S_2(c)$, proving
\[
1-d_{F,\mathcal{G}_1\square\mathcal{G}_2,\mathcal{P}}(\mathbf{g},\mathbf{h})=S_1(r)S_2(c).
\]
The equivalent formula for $d$ follows from
\[
1-S_1S_2=(1-S_1)+(1-S_2)-(1-S_1)(1-S_2).
\]
\end{proof}

\newpage 
\section{Some examples from numerical experiments}\label{S:Examples}

{
We illustrate the metric using synthetic graphs.
The histograms provide qualitative illustrations of how the graph structure affects the distribution of distances.
In the continuous-coordinate experiments, vectors are sampled uniformly from $[0,1]^N$, and the base distance is $d_i(g_i,h_i)=\frac12|g_i-h_i|$.
All present edges have weight $p_{ij}=1$.
Thus the distance used in these experiments is
\[
d_{\mathcal{X},\mathcal{G},\mathcal{P}}(\mathbf g,\mathbf h)
=1-\frac1N\sum_{j=1}^N\prod_{i\in B_j(\mathcal G)}
\left(1-\frac12|g_i-h_i|\right).
\]
The discrete poset specialization is discussed separately below.

}

\begin{figure*}[th]
  \centering
  \vspace{-0.1in}
\includegraphics[width=0.90\textwidth]{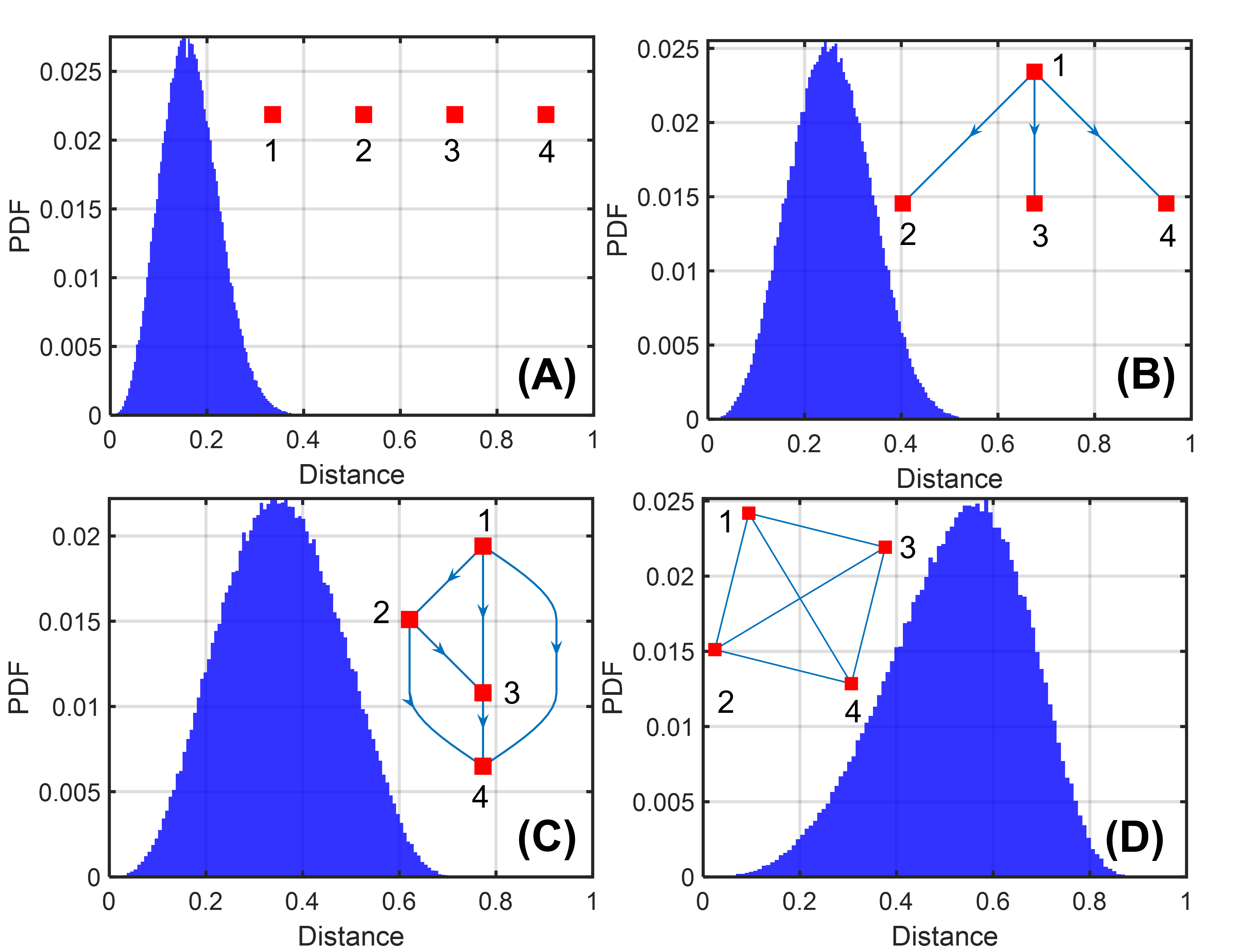}  
  \caption{Distribution of $d_{\mathcal{X},\mathcal{G},\mathcal{P}}$ for different graphs that model the dependencies between the elements of $\mathcal{X}$ generated by uniformly sampling from the volume of the cube $\left[0,1\right]^4$:(A) Fully disconnected graph (no dependency); (B) a star directed graph; (C) directed acyclic graph; and (D) an undirected graph.}  
  \label{fig:smallgraph}
\end{figure*}

{
Figure~\ref{fig:smallgraph} compares four graphs on $N=4$ vertices.
The distance distribution shifts to the right as edges are added, in agreement with Corollary~\ref{intro:C01}.
The different shapes illustrate the dependence of the distribution on the graph. They do not by themselves identify a unique topology.
}

\begin{figure*}[th]
  \centering
  \vspace{-0.1in}
\includegraphics[width=0.75\textwidth]{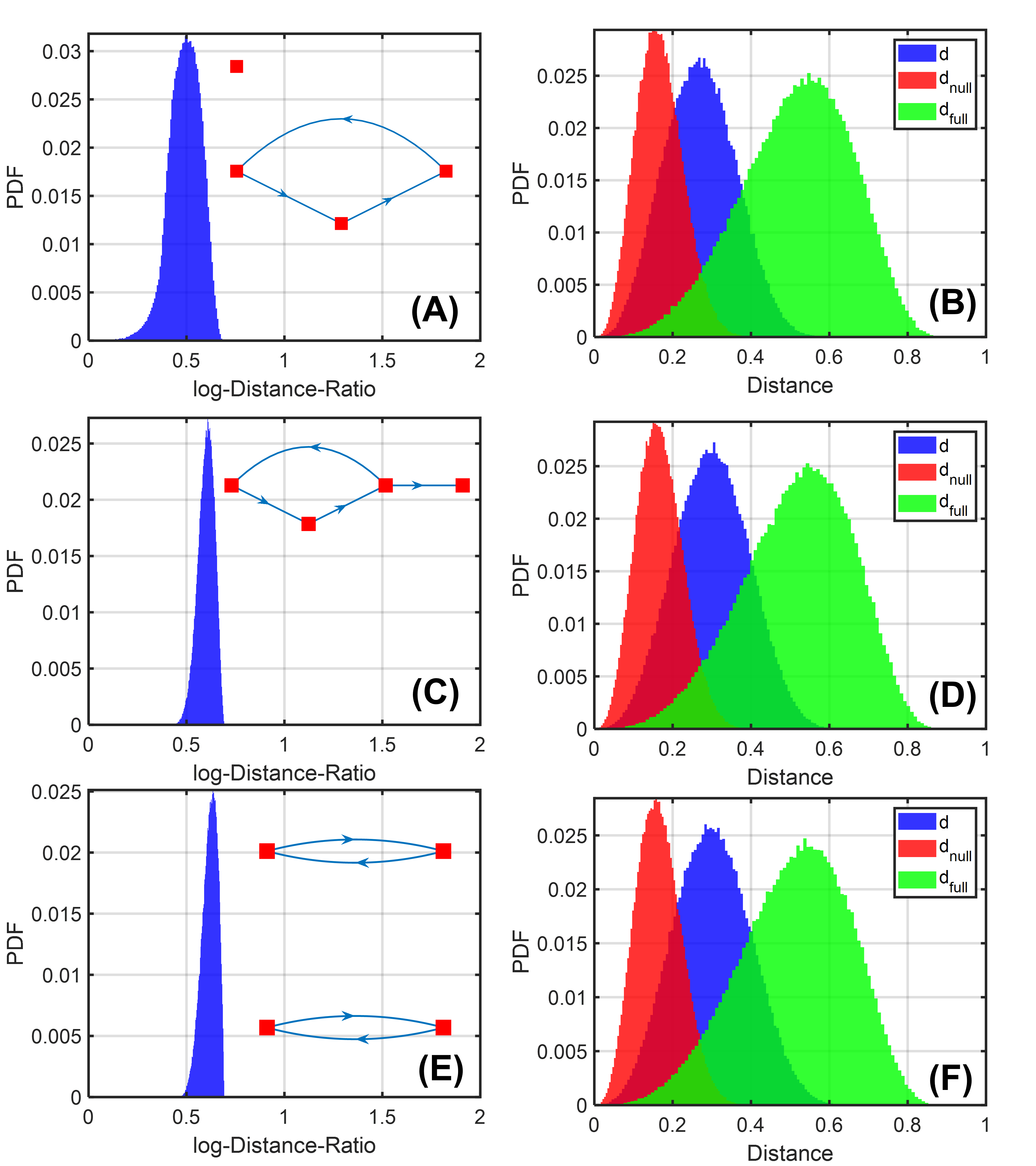}  
  \caption{Distribution of log-distance-ratio of $d_{\mathcal{X},\mathcal{G},\mathcal{P}}$ with respect to the $d_{null}$ and $d_{full}$ for different graphs shown in the inset.}
  \label{fig:incrementalgraph}
\end{figure*}

{
Figure~\ref{fig:incrementalgraph} compares distances as edges are added to a graph on four vertices containing a directed cycle.
We use the edgeless and complete bidirected graphs as references, with distances $d_{null}$ and $d_{full}$.
For every fixed pair of vectors, Corollary~\ref{intro:C01} gives $d_{null}\leq d\leq d_{full}$.
For pairs with $d_{null}>0$, the log-distance-ratio $\log(d/d_{null})$ is therefore nonnegative.
}

{
The ratio compares the graph distance with the edgeless distance for the same pair of vectors, reducing the influence of the pair's overall coordinate separation.
The logarithm expresses this relative increase on an additive scale.
This can make differences between the plotted graph distributions easier to see than raw distances, especially when their overall scales are similar.
It is a descriptive comparison, not a guarantee of better separation for every pair of graphs.
Pairs with $d_{null}=0$ are excluded from this ratio, and very small denominators require care in numerical evaluation.
}

\begin{figure*}[th]
  \centering
  \vspace{-0.1in}
\includegraphics[width=0.75\textwidth]{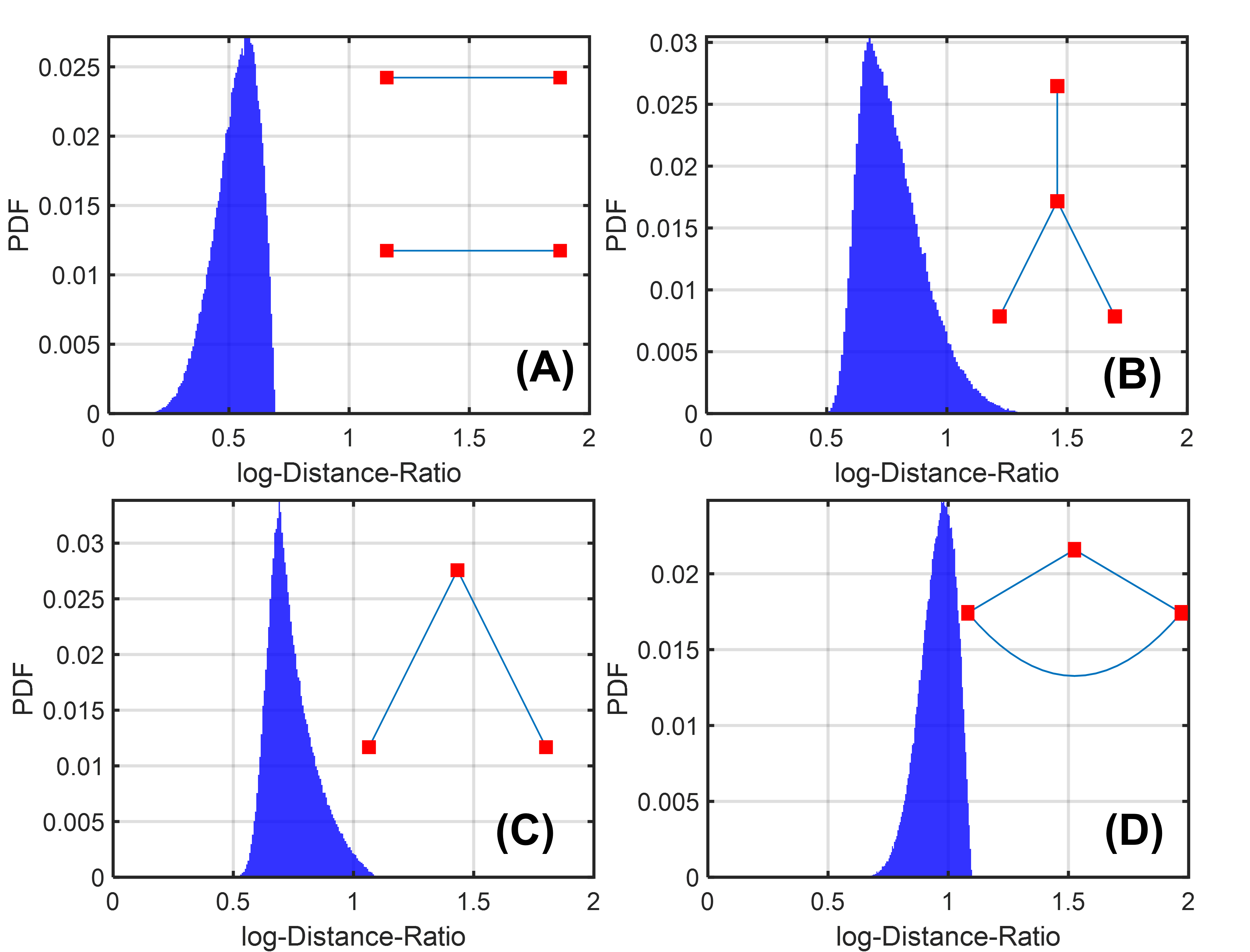}  
  \caption{Distribution of log-distance-ratio of $d_{\mathcal{X},\mathcal{G},\mathcal{P}}$ with respect to the $d_{null}$ and $d_{full}$ for different undirected graphs with different number of vertices (shown in the inset).}
  \label{fig:verticesgraph}
\end{figure*}

{
Figure~\ref{fig:verticesgraph} compares graphs with three and four vertices.
Panels (B) and (C) have more similar log-distance-ratio distributions than panels (A) and (D), despite their different numbers of vertices.
The comparison illustrates that the vertex count alone does not determine these distributions.
}

\begin{figure*}[thp]
  \centering
  \vspace{-0.1in}
\includegraphics[width=\textwidth]{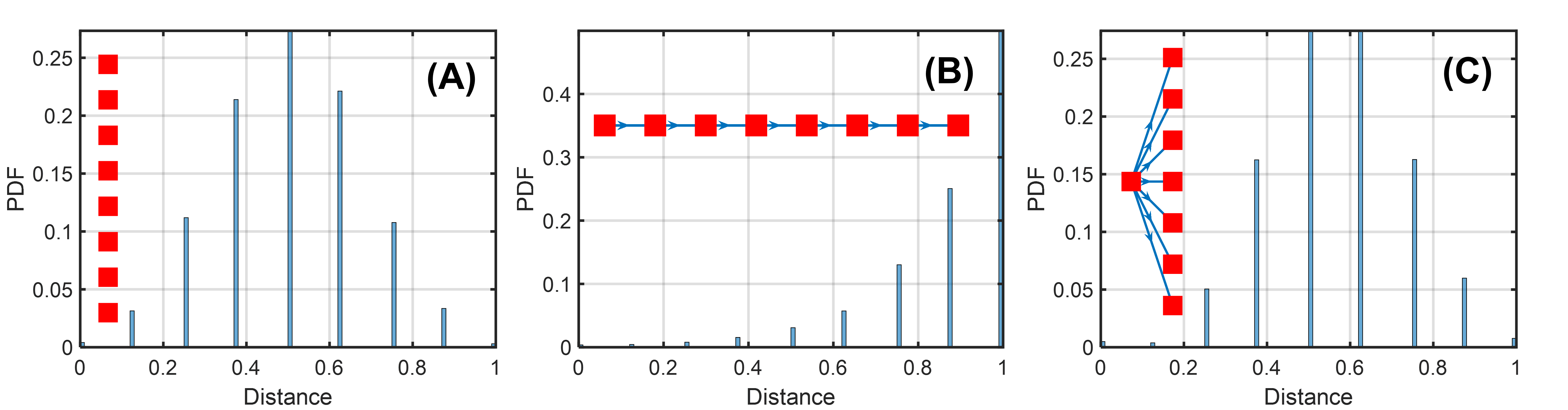}  
  \caption{Distribution of $d_{\mathcal{X},\mathcal{G},\mathcal{P}}$ corresponding to graphs with Hasse diagrams for (A) {an antichain}; (B) a Chain poset; and (C) a Star poset.}
  \label{fig:poset}
\end{figure*}

{
For the poset specialization in Corollary~\ref{intro:C1}, the sampling space is the set of vertices $\{0,1\}^N$ of the unit cube, and the base distance is the discrete metric $d_i(g_i,h_i)=|g_i-h_i|$.
The graph used in the metric is $G^{\preceq}$, with the reachability convention of that corollary.
These conventions give the normalized poset distance $|\langle\operatorname{supp}(g-h)\rangle|/N$.
Figure~\ref{fig:poset} is intended to illustrate the antichain, chain, and star cases.
In the antichain case, uniform independent binary vectors give the binomial distribution of normalized Hamming distance on the full binary space, rather than the weight distribution of a particular Hamming code.
The factor $1/2$ used in the continuous experiments does not give this poset specialization.
}

\begin{figure*}[thp]
  \centering
  \vspace{-0.1in}
\includegraphics[width=\textwidth]{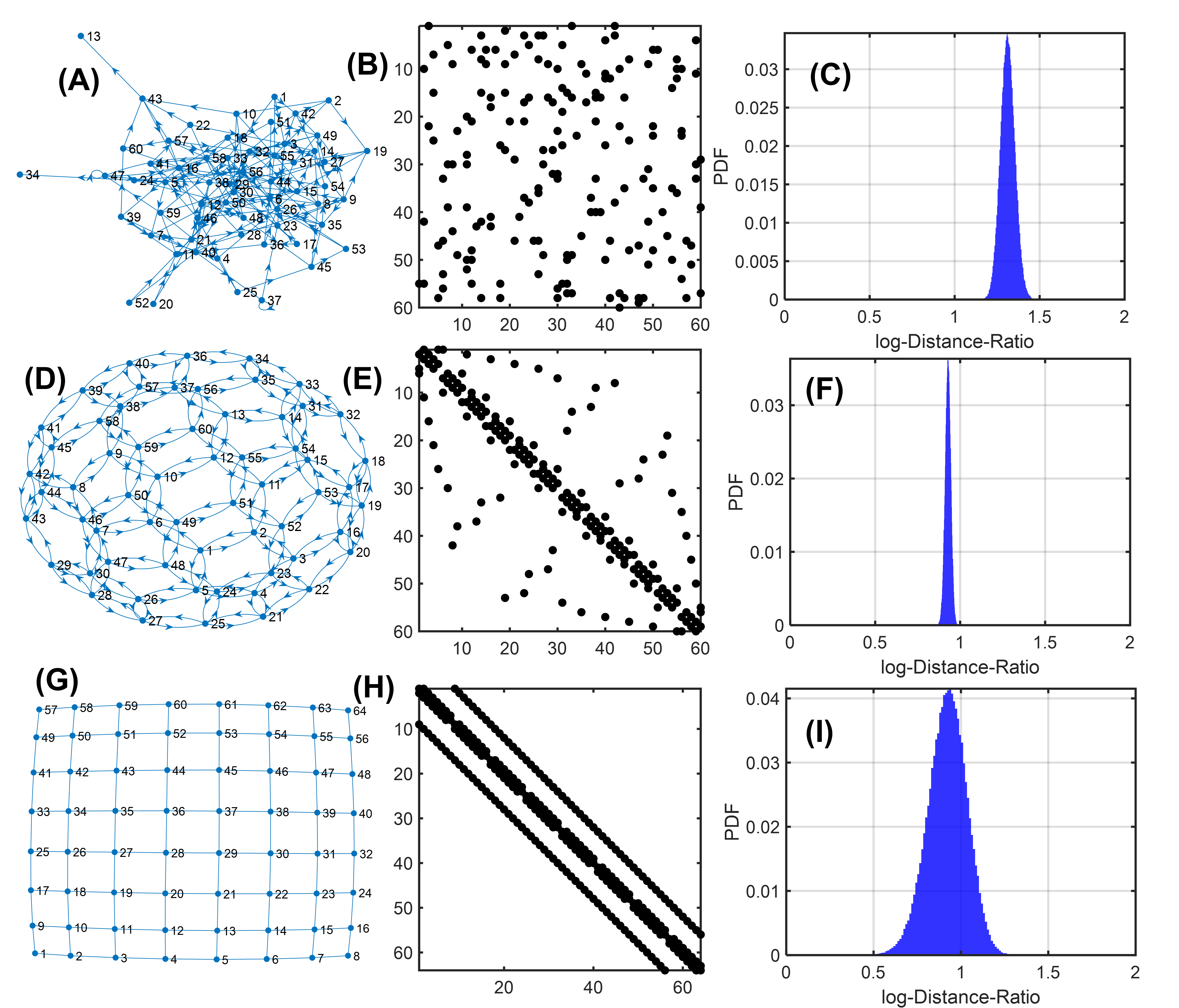}  
  \caption{Distribution of log-distance-ratio for different bidirectional graphs and their respective adjacency matrices: (A)-(C) random sparse graph; (D)-(F) graph with vertices arranged and connected over a 3-D bucky-ball; (G)-(I) graph with vertices arranged in a 2-D grid.}
  \label{fig:diffgraphs}
\end{figure*}

{
Figure~\ref{fig:diffgraphs} compares a random sparse graph, a bucky-ball graph, and a rectangular grid, together with their adjacency matrices.
The log-distance-ratio distributions differ across these examples, with the bucky-ball distribution more concentrated about its mean.
This concentration describes the spread of the relative distances. It does not by itself imply that the graph distance is close to the edgeless distance.
}

\begin{figure*}[thp]
  \centering
  \vspace{-0.1in}
\includegraphics[width=\textwidth]{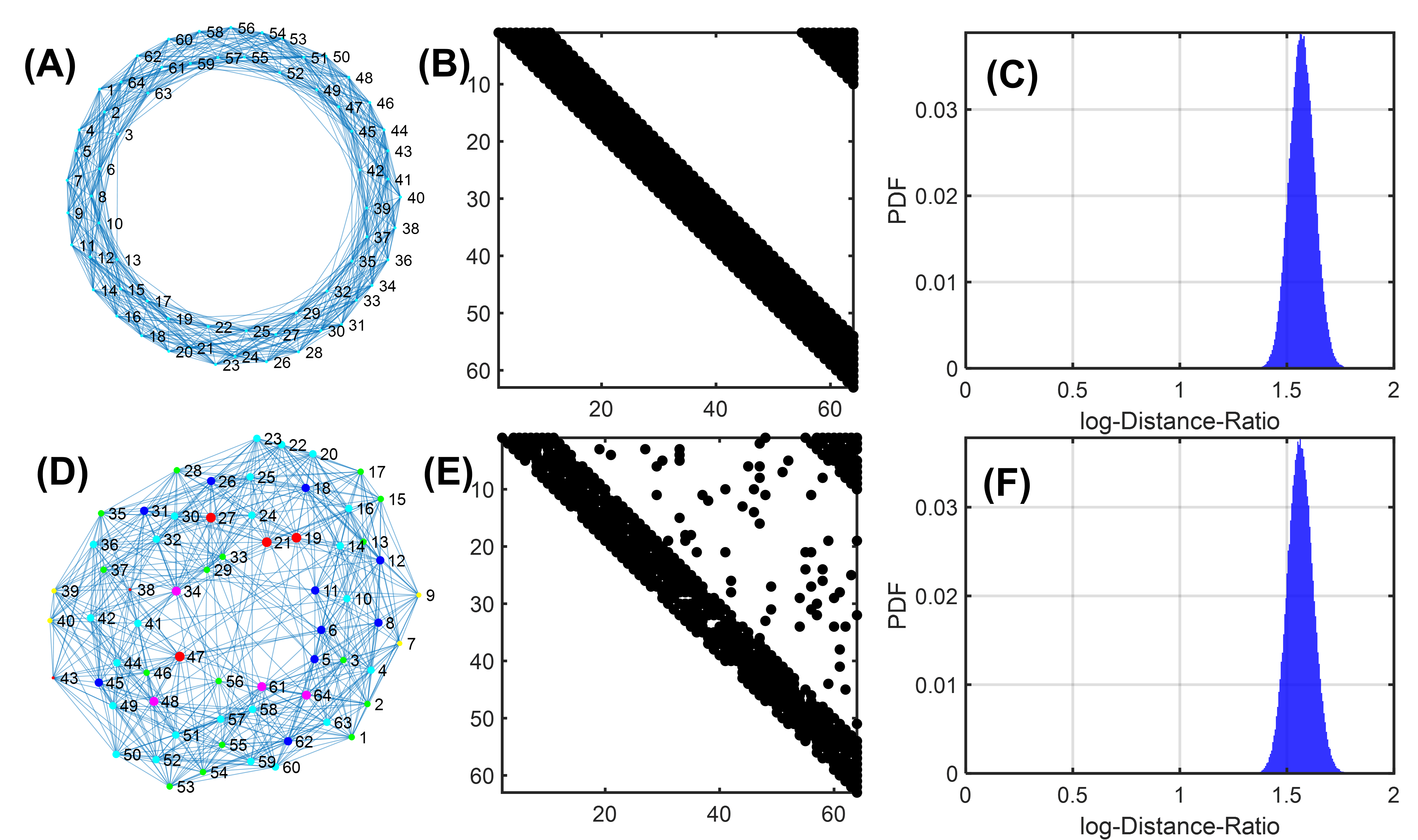}  
  \caption{Log-distance-ratio distribution for {Watts--Strogatz} graphs with 64 vertices and node degree = 10: (A) graph with $\beta=0$, (B) its adjacency matrix, and (C) log-distance-ratio distribution, (D) graph with $\beta = 0.025$, (E) its adjacency matrix, and (F) log-distance-ratio distribution. }
  \label{fig:strogatzconnectivity}
\end{figure*}

{
We next use the Watts--Strogatz model~\cite{Watts1998}, whose rewiring parameter $\beta$ controls departures from a regular lattice.
Figure~\ref{fig:strogatzconnectivity} compares the cases $\beta=0$ and $\beta=0.025$ with the same number of vertices and initial degree.
The similar log-distance-ratio distributions show that this comparison is not strongly affected by the small change in $\beta$ in these examples. They do not establish a general measure of clustering or small-world structure.
}

\begin{figure*}[thp]
  \centering
  \vspace{-0.1in}
\includegraphics[width=\textwidth]{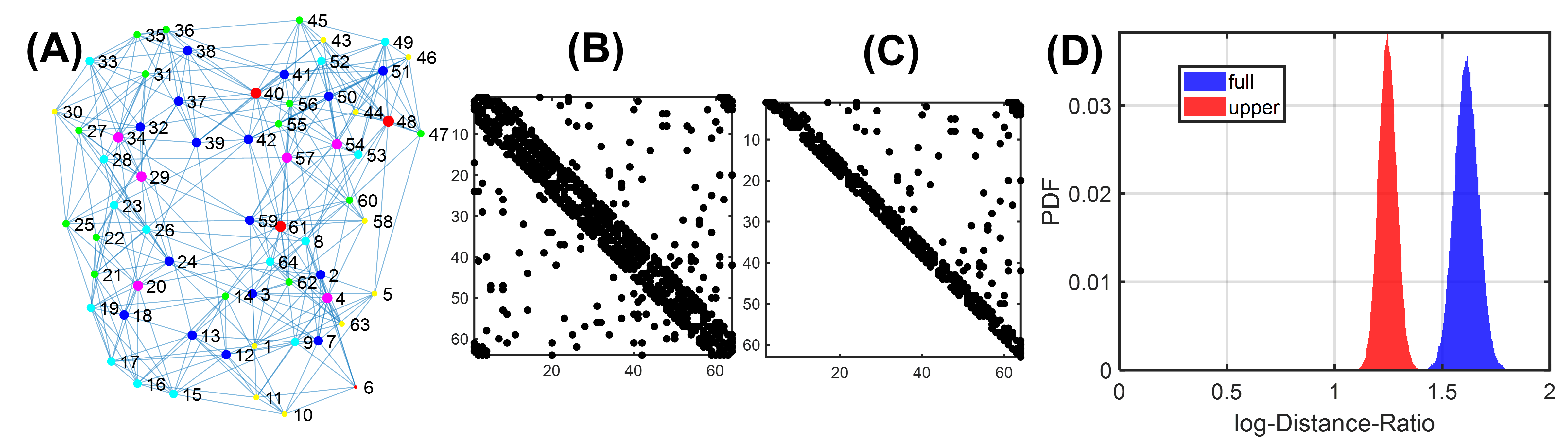}  
  \caption{Comparison of directed and undirected graph corresponding to a (A) Watts-Strogatz random graph with degree = 5, and $\beta = 0.2$. (B) adjacency matrix of the undirected graph; (C) adjacency matrix of the directed graph;(D) distribution of log-distance-ratios for the undirected graph (labeled as full) and {the directed graph} (labeled as upper).}
  \label{fig:strogatzdirection}
\end{figure*}

{
Figure~\ref{fig:strogatzdirection} compares the undirected graph with a directed version obtained from the upper triangular part of its adjacency matrix.
The directed graph has fewer edges, so Corollary~\ref{intro:C01} explains the lower mean log-distance-ratio observed in panel (D).
Panels (B) and (C) show the undirected and directed adjacency matrices, respectively.
}

\section{Expectations and mean distances}\label{S:EEM}

{
The purpose of this section is to give a mathematical basis for our numerical observations.
} 
We begin with a brief discussion of mean distances and how they shift under edge densification.
This clarifies the behavior observed in our empirical histograms and provides design
principles: (i) adding edges monotonically increases expected distance; (ii) under
coordinate-pair independence the mean factorizes into simple node-wise summaries; (iii) the effect
of any single new edge is uniformly bounded by $1/N$, giving robust, interpretable
comparative statics.

\medskip

\noindent\textbf{Expectation conventions.}
For a probability measure $\mu$ on $X\times X$ and any measurable $F:X\times X\to\mathbb{R}$,
\begin{equation*}
  \E_{(g,h)\sim\mu}\big[F(g,h)\big]
  \;:=\;\int_{X\times X} F(g,h)\,d\mu(g,h).
\end{equation*}
If $C\subset X$ is a finite code and $\mu$ is the uniform law on $C\times C$, then
\begin{equation*}
  \E_{(g,h)\sim\mu}\big[\dG(g,h)\big]
  \;=\;\frac{1}{|C|^2}\sum_{g,h\in C}\dG(g,h).
\end{equation*}

\medskip

\noindent\textbf{Standing notation.}
Let $G$ be a directed graph on vertices $[N]=\{1,\dots,N\}$ with positive edge weights
$P=(p_{ij})_{(i,j)\in E(G)}\subset(0,1]$. For $j\in[N]$ write $\Bin{j}:=\{j\}\cup\{i:(i,j)\in E(G)\}$ and set $p_{jj}=1$.
For metric spaces $(X_i,d_i)$ with $d_i\colon X_i\times X_i\to[0,1]$ and points
$g=(g_1,\dots,g_N)$, $h=(h_1,\dots,h_N)$ in $X:=\prod_i X_i$, set
\begin{equation}\label{eq:graph-metric}
  \dG(g,h)
  \;:=\;1-\frac{1}{N}\sum_{j=1}^{N}\;\prod_{i\in \Bin{j}} \bigl(1-d_i(g_i,h_i)\bigr)^{\!1/p_{ij}}\,.
\end{equation}
Let $\mu$ be any probability measure on $X\times X$ and define the mean distance
\begin{equation}\label{eq:mean-def}
  \mathbf{M}_{\mu}(G,P)\;:=\;\E_{(g,h)\sim\mu}[\,\dG(g,h)\,]\,.
\end{equation}

\begin{Theorem}[Monotonicity under edge densification]\label{thm:mean-monotone}
Fix $G$ and $P$ and consider $G'$ obtained from $G$ by adding any collection of new edges
(on the same vertex set), assigning weights in $(0,1]$ to the new edges while keeping the old
weights unchanged. Then for every $(g,h)$ one has $\dG(g,h)\le d_{X,G',P}(g,h)$, and consequently
\begin{equation*}
  \mathbf{M}_{\mu}(G,P)\;\le\;\mathbf{M}_{\mu}(G',P)\,.
\end{equation*}
Moreover, if a single edge $(i\to j)$ with weight $p\in(0,1]$ is added, then with
$Q_j(G):=\prod_{k\in \Bin{j}}\bigl(1-d_k(g_k,h_k)\bigr)^{\!1/p_{kj}}$ one has the exact
increment identity
\begin{equation}\label{eq:delta-mean-single-edge}
  d_{X,G',P}(g,h)-\dG(g,h)
  \;=\;\frac{1}{N}\,Q_j(G)\,\Bigl(1-\bigl(1-d_i(g_i,h_i)\bigr)^{\!1/p}\Bigr)\;\ge 0\,.
\end{equation}
Averaging \eqref{eq:delta-mean-single-edge} over $\mu$ yields
$\mathbf{M}_{\mu}(G',P)-\mathbf{M}_{\mu}(G,P)=\E_{\mu}[\,\text{RHS}\,]$.
\end{Theorem}
\begin{proof}
Fix $(g,h)$. For each $j$, set
\[
S_j(G)\;:=\;\prod_{k\in\Bin{j}}\bigl(1-d_k(g_k,h_k)\bigr)^{1/p_{kj}}\in[0,1].
\]
If a single edge $(i\to j)$ with weight $p\in(0,1]$ is added, then
\[
S_j(G') \;=\; S_j(G)\cdot \bigl(1-d_i(g_i,h_i)\bigr)^{1/p}\ \le\ S_j(G),
\]
since the extra factor lies in $[0,1]$. All other $S_\ell$ are unchanged.
Therefore, we have
\[
d_{X,G',P}(g,h)-d_{X,G,P}(g,h)
= \frac{1}{N}\big(S_j(G)-S_j(G')\big)
= \frac{1}{N}S_j(G)\Bigl(1-\bigl(1-d_i(g_i,h_i)\bigr)^{1/p}\Bigr)\ \ge 0,
\]
which is exactly \eqref{eq:delta-mean-single-edge}. Adding multiple edges one by one yields the inequality 
$d_{X,G',P}(g,h)\ge d_{X,G,P}(g,h)$ for all $(g,h)$. Finally, taking expectations gives
$\mathbf{M}_\mu(G,P)\le \mathbf{M}_\mu(G',P)$.
This completes the proof of our theorem. 
\end{proof}

The effect of changing weights on existing edges can be summarized as follows. 

\begin{Proposition}\label{P:weight-monotone}
For fixed $(g,h)$ and edge set, the function $p\mapsto 1-\bigl(1-d_i(g_i,h_i)\bigr)^{1/p}$ is decreasing on $(0,1]$.
{
Hence, if $0<\Delta<p_{ij}$ and $p_{ij}+\Delta\leq1$, then
} 
\[
  d_{X,G,P+\Delta}\, \le \, d_{X,G,P}\, \le\, d_{X,G,P-\Delta}.
\]
\end{Proposition}
\begin{proof}
Write $x:=1-d_i(g_i,h_i)\in[0,1]$.
{
For $x=0$ or $x=1$, the function $f(p):=x^{1/p}$ is constant. For $0<x<1$,
} we have
$f'(p)=x^{1/p}(-\ln x)/p^2\ge0$. Thus $f$ increases and $1-f$ decreases in $p$.
\end{proof}

We now describe the closed form of $\mathbf{M}_{\mu}(G,P)$ under coordinate-pair independence.

\begin{Theorem}\label{thm:mean-factorization}
Assume $\{(g_i,h_i)\}_{i=1}^N$ are independent under $\mu$. 
For $t>0$ set
\[
m_i(t):=\E\bigl[(1-d_i(g_i,h_i))^{t}\bigr]\in[0,1].
\]
Then we have 
\begin{equation}\label{eq:mean-closed-form}
  \mathbf{M}_{\mu}(G,P)
  \;=\;1-\frac1N\sum_{j=1}^N\;\prod_{i\in\Bin{j}} m_i\!\Bigl(\frac{1}{p_{ij}}\Bigr).
\end{equation}
If, in addition, the coordinate pairs $(g_i,h_i)$ are identically distributed, with common $m(t)$, and all edge weights are equal
to $p$, then
\begin{equation}\label{eq:mean-iid}
  \mathbf{M}_{\mu}(G,P)
  \;=\;1-\frac{1}{N}\sum_{j=1}^N m(1)\,m\!\Bigl(\frac{1}{p}\Bigr)^{\indeg(j)}\!,
  \qquad
  \Delta\mathbf{M}_{\mu}\;=\;\frac{1-m(1/p)}{N}\,m(1)\,m\!\Bigl(\frac{1}{p}\Bigr)^{\indeg(j)}
\end{equation}
for the single-edge increment at a node $j$.
\end{Theorem}
\begin{proof}
Independence gives $\E\bigl[\prod_{i\in\Bin{j}}(1-d_i)^{1/p_{ij}}\bigr]=\prod_{i\in\Bin{j}}m_i(1/p_{ij})$, which
inserted into \eqref{eq:graph-metric} yields \eqref{eq:mean-closed-form}. The i.i.d. case and increment
follow as stated.
\end{proof}


We close this section by presenting two interesting consequences of our Theorems~\ref{thm:mean-monotone} and~\ref{thm:mean-factorization}.

\begin{Corollary}[Edge sensitivity and a telescoping bound]\label{cor:edge-sensitivity}
For any single added edge, the pointwise increment in \eqref{eq:delta-mean-single-edge} is bounded by $1/N$.
Consequently, for any two graphs $G\subseteq G'$ on the same vertex set,
\[
  0\;\le\;\mathbf{M}_{\mu}(G',P)-\mathbf{M}_{\mu}(G,P)\;\le\;\frac{|E(G')\setminus E(G)|}{N}\,.
\]
\end{Corollary}
\begin{proof}
In \eqref{eq:delta-mean-single-edge}, $Q_j(G)\le1$ and
$0\le 1-(1-d_i)^{1/p}\le1$, hence the $1/N$ bound. Sum over added edges.
\end{proof}

\begin{Corollary}[Jensen-type bounds for the mean distance]\label{cor:mean-bounds}
For $t\ge1$ the map $x\mapsto x^{t}$ is convex and satisfies $x^{t}\le x$ on $[0,1]$. Thus
\[
  \prod_{i\in \Bin{j}}\Big(\E[1-d_i]\Big)^{\!1/p_{ij}}
  \;\le\;
  \prod_{i\in \Bin{j}} m_i\!\Big(\tfrac{1}{p_{ij}}\Big)
  \;\le\;
  \prod_{i\in \Bin{j}}\E[1-d_i].
\]
{Under the independence assumption of Theorem~\ref{thm:mean-factorization}, substitution into \eqref{eq:mean-closed-form} gives bounds on $\mathbf{M}_{\mu}(G,P)$. Their directions are reversed because the mean distance is one minus the averaged product.}
\end{Corollary}
\begin{proof}
Fix $i$ and set $X_i:=1-d_i(g_i,h_i)\in[0,1]$. Since $p_{ij}\in(0,1]$, we have $t_{ij}:=1/p_{ij}\ge1$.
The map $x\mapsto x^{t_{ij}}$ is convex on $[0,1]$ (since $t_{ij}(t_{ij}-1)x^{t_{ij}-2}\ge0$), hence by Jensen's inequality
\[
  \big(\E[X_i]\big)^{t_{ij}} \;\le\; \E\!\big[X_i^{\,t_{ij}}\big].
\]
Also, for $x\in[0,1]$ and $t_{ij}\ge1$ we have $x^{t_{ij}}\le x$, so $\E[X_i^{\,t_{ij}}]\le \E[X_i]$.
But it is easy to check that 
\[
  \big(\E[1-d_i]\big)^{1/p_{ij}}
  \;\le\;
  m_i\!\Big(\tfrac{1}{p_{ij}}\Big)
  \;\le\;
  \E[1-d_i].
\]
Then multiplying over $i\in \Bin{j}$ gives the claimed product inequalities, and inserting into
\eqref{eq:mean-closed-form} yields the bounds on $\mathbf{M}_{\mu}(G,P)$. 
This completes the proof. 
\end{proof}

\begin{Example}\label{ex:binary}
Suppose that $d_i(g_i,h_i)\in\{0,1\}$. For example, $d_i$ can be the discrete metric $\mathbf{1}_{\{g_i\ne h_i\}}$. 
Then $(1-d_i)^t=\mathbf{1}_{\{g_i=h_i\}}$ for all $t>0$, so $m_i(t)=P[g_i=h_i]$ is independent of $t$. If the coordinates are independent with match probabilities $\rho_i$, then we have 
\[
  \mathbf{M}_{\mu}(G,P)
  \;=\;1-\frac{1}{N}\sum_{j=1}^{N}\prod_{i\in \Bin{j}}\rho_i.
\]
If, in addition, the coordinate pairs are i.i.d. with $\rho_i\equiv\rho$, then the formula for $\mathbf{M}_{\mu}(G,P)$ becomes
$\mathbf{M}_{\mu}(G,P)=1-\frac{1}{N}\sum_{j=1}^{N}\rho^{\,\indeg(j)+1}$.
This formulation makes the impact of in-degrees on the mean completely explicit.
\end{Example}

\section{Conclusions}\label{S:Conclusions}

In this paper, we introduced a novel metric space for measuring distances using arbitrary directed graphs, making it applicable to graphs and sub-graphs representing diverse discrete objects. Using synthetic graphs, we demonstrated how statistical measures derived from our proposed distance metrics can illuminate various graph properties.

We believe this metric space will be fundamental to numerous domains where distance-based algorithms are routinely employed for visualization and analysis, such as clustering algorithms in machine learning. Consequently, for certain problems and data structures, our framework could eliminate the need for developing new mathematical approaches, particularly when training data for AI models is sparse. This relevance can be highlighted by the success of the Wasserstein distance, which has been used to improve the performance of generative AI models~\cite{Lu2023}.

\section*{Acknowledgements.}
M.B.C gratefully acknowledges the support provided by the Louisiana Board of Regents award LEQSF(2023-25)-RD-A-21 for this research. S.C acknowledges the support from the National Science Foundation grants ECCS: 2332166, FET: 2208770 and CNS-FuSE2: 2425444.
The authors thank the referees for their careful reading of the article and constructive comments.

\bibliographystyle{plain}
\bibliography{references}
\end{document}